\documentclass[a4paper,12pt]{article}
\usepackage{amsfonts}
\usepackage{amsmath}
\usepackage{amssymb}
\usepackage{amsthm}
\usepackage{url}
\title{Non-hyperelliptic curves of genus three over finite fields of characteristic two}
\author{Enric Nart\thanks{
acknowledges financial support of the spanish MECD and MCYT
through the ``Servicio de Acciones de Promoci\'on y Movilidad" and
the project BFM-2003-06092} \ and Christophe Ritzenthaler
\thanks{acknowledges the financial support provided through the European
Community's Human Potential Programme under contract
HPRN-CT-2000-00114, GTEM }}

\date{}

\newtheorem{teor}{Theorem}[section]
\newtheorem{cor}[teor]{Corollary}

\newtheorem{lem}[teor]{Lemma}
\newtheorem{prop}[teor]{Proposition}

\newcommand{\C}{\mathbb C}
\def\op{\operatorname}

\def\aa{{\mathcal A}}

\def\al{\alpha}
\def\as{\op{AS}}
\def\aut{\op{Aut}}
\def\bb{{\mathcal B}}
\def\be{\bigskip}

\def\cc{{\mathcal C}}

\def\cq{{\mathcal C}_{\negmedspace{\mathcal Q}}}

\def\dd{{\mathcal D}}
\def\ddg{{\mathcal D}_{\ga}}

\def\dsc{\op{dsc}}

\def\diso{\lower.4ex\hbox{$\downarrow$}\raise.4ex\hbox{\mbox{\scriptsize $\wr$}}}
\def\e{\medskip}

\def\eps{\epsilon}

\def\fq{\mathbb F_q}
\def\ft{\mathbb F_2}
\def\g{\Gamma}
\def\ga{\gamma}

\def\gen#1{\big\langle\, {#1} \,\big\rangle}
\def\gga{\Gamma_{\gamma}}
\def\gl#1#2{\op{GL}_{#1}(#2)}

\def\imp{\,\Longrightarrow\,}

\def\iso{\,\stackrel{\mbox{\tiny $\sim\,$}}{\lra}\,}

\def\kb{\bar{k}}

\def\lg{l\raise.6ex\hbox to.2em{\hss.\hss}l}
\def\lra{\longrightarrow}

\def\md#1{\ \mbox{\rm(mod }{#1})}
\def\mm{{\mathcal M}}
\def\ni{\noindent}
\def\nn{{\mathcal N}}

\def\nv{N(v)}
\def\oo{{\mathcal O}}

\def\orb{\hbox to  .3em{$\backslash$}\backslash}

\def\pg#1{\op{PGL}_3(#1)}

\def\pr#1{\mathbb P^#1}
\def\qq{{\mathcal Q}}

\def\sg{\sigma}
\def\si{\,^{\sigma}\!}
\def\sij#1{\,^{\sigma^{#1}}\!}
\def\sii{\,\Longleftrightarrow\,}

\def\ss{{\mathcal S}}

\def\t{\,^t\!}

\def\tq{\,\,|\,\,}
\def\tr{\op{Tr}}

\def\uu{{\mathcal U}}

\def\ww{{\mathcal W}}

\def\z{\zeta}

\newcounter{cs}
\stepcounter{cs}
\newcommand{\casos}{\begin{itemize}}
\newcommand{\fcasos}{\end{itemize}\setcounter{cs}{1}}
\newcommand{\cas}{\item[(\alph{cs})]\stepcounter{cs}}
\newfont{\tit}{cmr12 scaled \magstep3}

\begin{document}

\maketitle

\begin{abstract}
Let $k$ be a finite field of even characteristic. We obtain in
this paper a complete classification, up to $k$-isomorphism, of
non singular quartic plane curves defined over $k$. We find
explicit rational normal models and we give closed formulas for
the total number of $k$-isomorphism classes. We deduce from these
computations the number of  $k$-rational points of the different
strata by the Newton polygon of the non hyperelliptic locus
$\mm_3^{\textrm{nh}}$ of the moduli space $\mathcal{M}_3$ of
curves of genus $3$. By adding to these computations the results
of \cite{ns}, \cite{scholten} on
 the hyperelliptic locus  we obtain a complete picture of these
 strata for $\mm_3$.
\end{abstract}

\section*{Introduction}
The study of non singular plane quartics has a long history, going
back to XIX-th century geometers, and it is still today a fruitful
area of research. These genus $3$ curves are the first example of
non hyperelliptic curves and their geometry is completely
different: classification and invariants of hyperelliptic curves
rely on a binary theory given by abscissas of Weierstrass points;
in the case of quartics, we deal with bitangents and we obtain a
ternary theory. Over $\C$, in spite of recent progress, the
des\-cription of the moduli space $\mm_3$ of genus $3$ curves is
not complete: from a geometric point of view, the works of
 L. Caporaso and E. Sernesi \cite{caporaso} and D. Lehavi \cite{lehavi}
show that the old construction of Riemann of a quartic from one of
its Aronhold systems (certain subsets of seven bitangents) give in
fact an explicit finite morphism between an open subset of
$\pr{6}$ and $\mm_3^{\textrm{nh}}$, the non hyperelliptic locus of
$\mm_3$. From the point of view of the theory of invariants, one
tries to understand the graded algebra $A$  of invariants for the
natural action of $\textrm{SL}(3,\C)$ on the vector space of
homogeneous polynomials of degree $4$
 in three  variables. Dixmier in \cite{dixmier} constructed
 a homogeneous system of parameters for $A$, of degrees $3, 6, 9, 12, 15, 18$ and
 $27$, but the precise number of generators of $A$ is not known.
At last, Katsylo in \cite{katsylo} proved that $\mathcal{M}_3$ is
rational.

Can we precise the situation over other fields? If $k$ is a field
of characteristic $2$, we may consider a stratification of $\mm_3$
by the $2$-rank of the Jacobian. One can show that quartics with
$2$-rank $3$ (resp. $2,1,0$) have $7$ (resp. $4,2,1$) bitangents.
For each case, C.T.C Wall \cite{wa} found a family of plane
quartics containing all $\kb$-isomorphy classes, together with a
group acting fully by isomorphisms on the family.

In this paper we obtain, for $k=\fq$ a finite field of
characteristic $2$, a complete classification up to
$k$-isomorphism of non singular quartics defined over $k$. Note
that, even for $k=\ft$, this had not been completely fulfilled
(cf. \cite{du} and the references quoted there). By applying
descent theory to each of the  four Wall families, we obtain for
each value of $n=7,4,2,1$ an explicit description of the set of
$k$-isomorphy classes of curves with $n$ bitangents in the form
$\amalg_i (G_i \backslash \nn_i^{(n)})$, where $\nn_i^{(n)}$ are
different families of rational normal models and each $G_i$ is a
finite group acting fully by $k$-isomorphisms on $\nn_i^{(n)}$.
Actually, we use the notation $\oo$ (for ordinary) as an
alternative to $\nn^{(7)}$ and $\ss$ (for supersingular) as an
alternative to $\mathcal{N}_0^{(1)}$. Altogether, we obtain
thirteen families of rational normal models, which are denoted by
$$\oo_1,\, \oo_2,\,\oo_3,\,\oo_4,\,\oo_{7,0},\,\oo_{7,1};\,\nn_1^{(4)},\,\nn_2^{(4)},\,
\nn_3^{(4)};\,\nn_1^{(2)},\,\nn_0^{(2)};\,\nn_1^{(1)},\,\ss.$$ We
compute also the $k$-automorphism group of each curve $C$ in
$\nn_i^{(n)}$, which is given by the stabilizer of $C$ under the
action of $G_i$. Moreover, we obtain closed formulas for the
number $|G_i \backslash \nn_i^{(n)}|$ of $k$-isomorphy classes in
each family $\nn_i^{(n)}$.

This is carried out in Section \ref{sec:1} for the ordinary case,
which serves as a prototype, and in Section \ref{sec:2} for the
non-ordinary cases. In Section \ref{sec:1} we precise also the
structure of the ordinary stratum of $\mm_3^{\textrm{nh}}$. We
show that it is isomorphic to an explicit open set of an affine
variety whose coordinate ring $R$ is the ring of invariants of
conics under the natural action of $\pg{\ft}$ on $\pr{2}(\kb)$. A
complete description of $R$ is obtained in \cite{mr}.

Finally, in Section \ref{sec:3} we deduce from these computations
the number of $k$-rational points of the strata by the Newton
polygon of $\mm_3^{\textrm{nh}}$. In all cases the Newton polygon
is determined by the 2-rank, except for the curves with 2-rank
zero, whose Newton polygon has either two sides with slopes 1/3,
2/3 (type 1/3) or one side with slope 1/2 (supersingular case). By
adding to these computations the results of \cite{ns},
\cite{scholten} on
 the hyperelliptic locus $\mm_3^{\textrm{h}}$ we obtain a complete
 picture of these strata for $\mm_3$. The results are:

\renewcommand\arraystretch{1.4}
\centerline{\begin{tabular}{c|ccccc}
&ordinary& 2-rank two&2-rank one&type $1/3$&supersingular\\
\hline$\mm_3^{\textrm{nh}}$& $q^6-q^5+1$&$q^5-q^4$&$q^4-q^3$&$q^3-q^2$&$q^2$\\
$\mm_3^{\textrm{h}}$&$q^5-q^4$&$q^4-2q^3+q^2$&$2(q^3-q^2)$&$q^2$&0\\
$\mm_3$& $q^6-q^4+1$&$q^5-2q^3+q^2$&$q^4+q^3-2q^2$&$q^3$&$q^2$
\end{tabular}}

\section{Ordinary curves}\label{sec:1}

Let $k=\fq$ be a finite field of characteristic 2 and let $\kb$ be
a fixed algebraic closure of $k$. We denote by $\qq$ the set of
quadratic forms:
\begin{equation}\label{qf}
Q(x,y,z)=ax^2+by^2+cz^2+dxy+eyz+fzx,
\end{equation}
with coefficients in $\kb$, such that:
$$
abc\ne 0,\ a+b+d\ne 0,\ b+c+e\ne 0, \ a+c+f\ne 0, \
a+b+c+d+e+f\ne1.
$$

To any $Q\in\qq$ we can associate the non-singular quartic plane
curve $C_Q$ determined by the equation:$$ C_Q\colon\quad
Q(x,y,z)^2=xyz(x+y+z)\qquad (\mbox{ shortly written as
}Q^2=xyz(x+y+z))
$$

All these curves have the same set of bitangents. It is the Fano
plane $\bb_0$ formed by the seven lines of $\pr{2}$ defined over
$\ft$. If we identify the lines of $\pr{2}$ with linear forms in
$x,y,z$ up to multiplication by a non-zero constant, we can write:
$$
\bb_0:=\{x,\ y, \ z, \ x+y+z, \ x+y, \ y+z, \ x+z\}.
$$

Let us denote by $\cq$ the family of all these curves $C_Q$. This
family contains all $\kb$-isomorphy classes of ordinary
non-hyperelliptic curves of genus three \cite{wa}, \cite{r}.

\begin{prop}\label{pr1}
Let $C$ be a non-singular quartic plane curve defined over $\kb$.
The following conditions are equivalent:
\begin{enumerate}
\item The Jacobian variety $J_C$ of $C$ is ordinary.
\item $C$ has seven bitangents.
\item The set $\bb$ of bitangents of $C$ is a Fano plane; that is,
\begin{equation}\label{fano}
\bb=\{\ell_1,\,\ell_2,\,\ell_3,\,\ell_1+\ell_2+\ell_3,\,\ell_1+\ell_2,\,\ell_2+\ell_3,\,
\ell_1+\ell_3\},
\end{equation}for some linear forms $\ell_1,\ell_2,\ell_3$.
\item The set of bitangents of $C$ is $\pg{\kb}$-equivalent to the Fano plane $\bb_0$.
\item $C$ is isomorphic to some curve $C_Q$, with $Q\in\qq$.\qed
\end{enumerate}
\end{prop}

For any isomorphism $\phi\colon C\to C_Q$, given by three linear
forms
$$\phi(x,y,z)=(\ell_1(x,y,z),\ell_2(x,y,z),\ell_3(x,y,z)),$$
the set of bitangents of $C$ is the Fano plane generated by
$\ell_1$, $\ell_2$, $\ell_3$ as in (\ref{fano}). Hence,
\begin{lem}\label{degdd}
Let $k\subseteq K\subseteq \kb$ be a finite extension of $k$. An
ordinary curve $C$ has all its bitangents defined over $K$ iff
 $C$ is $K$-isomorphic to $C_Q$ for some $Q\in\qq$ with coefficients in $K$.\qed
\end{lem}

If $C$ is defined over $k$, the minimum field of definition of all
bitangents of $C$ has degree 1,2,3,4 or 7 over $k$ \cite{r}.

In this section we describe the $k$-isomorphy classes of ordinary
non-singular quartics by applying descent theory to the family
$\cq$ (paragraph \ref{ddata}). We find explicit formulas for the
number of curves and we exhibit rational normal models for them
(paragraphs \ref{nbr} and \ref{ratmod}). Before of that, we recall
some generalities, we study the action of the subgroup of
$\pg{\kb}$ that preserves the family $\cq$ and we describe the
ring of invariants for ordinary curves (paragraph \ref{action}).

\subsection{Action of $\pg{\ft}$ on the curves $C_Q$ and invariants for ordinary curves}\label{action}
We have a natural action of $\pg{\kb}$ on the right on the set of
homogeneous polynomials $F(x,y,z)$ in three variables, up to
multiplication by a non-zero
 constant:
$$
F^{\ga}(x,y,z)=F(\ell_1(x,y,z),\ell_2(x,y,z),\ell_3(x,y,z)),\qquad
\forall \ga\in\pg{\kb},
$$where  $\ell_1,\ell_2,\ell_3$ are the linear forms having as coefficients the entries of the rows of
any representative of $\ga$ in $\gl{3}{\kb}$. If we think $\ga$ as
an automorphism of $\pr{2}$ and $V(F)\subseteq \pr{2}$ is the
subvariety of zeros of $F$, we have
$\ga(V(F))=V(F^{\ga^{-1}})$.\be

\ni{\bf Definition. }Let $\nn$ be a family of smooth projective
curves defined over a field $K\subseteq\kb$ and let $G$ be a group
acting on the left on $\nn$. We say that $G$ {\it acts by
$K$-isomorphisms} if there are distinguished $K$-isomorphisms
$g_C\colon C\iso g(C)$  for all pairs $C\in\nn$, $g\in G$,
satisfying:
$$
(g'g)_C=g'_{g(C)}\circ g_C,\qquad \forall C\in\nn,\,\forall
g,g'\in G.
$$We say that $G$ {\it acts fully by $K$-isomorphisms} if moreover:
$$
\op{Isom}_K(C,C')=\{g_C\tq g\in G,\, g(C)=C'\},\qquad \forall
C,C'\in\nn.
$$  In particular, in this latter case $\aut_K(C)$ can be identified with the stabilizer of $C$
under the action of $G$. In the case $K=\kb$ we say simply that
$G$ {\it acts (fully) by isomorphisms}. \be

If $C$ is a non-singular plane quartic, the inclusion $C\subseteq
\pr{2}$ coincides with the canonical embedding; hence, the group
$\pg{\kb}=\aut(\pr{2})$ acts fully by isomorphisms on the family
of all non-singular quartics, with $g_C=g_{\vert C}$ for all
$g\in\pg{\kb}$.

Let us denote by $\g:=\aut_{\ft}({\pr{2}})=\pg{\ft}= \gl{3}{\ft}$
the subgroup of $\pg{\kb}$ of those automorphisms that leave the
Fano plane $\bb_0$ invariant as a set. The group $\g$ acts fully
by isomorphisms on the family $\cq$. In fact, if we apply to the
curve $C_Q$ an automorphism $\ga\in\g$ we obtain $\ga(C_Q)=C_{Q'}$
for some $Q'\in\qq$. On the other hand, any isomorphism $C_Q\iso
C_{Q'}$ between two curves in $\cq$ sends the set of bitangents of
$C$ into the set of bitangents of $C'$; hence, it is given by an
element in $\g$.

Let us explicit the action on the left of $\g$ on $\qq$, that
reflects the natural action on the curves $C_Q$:

$$
\ga(Q):=Q^{\ga^{-1}}+H_{\ga^{-1}},\quad\forall \ga\in\g, \ \forall
Q\in\qq,
$$
where, for any $\ga$ with rows $\ell_1,\ell_2,\ell_3$, $H_{\ga}$
is the quadratic form determined by:
$$
\ell_1\ell_2\ell_3(\ell_1+\ell_2+\ell_3)=xyz(x+y+z)+H_{\ga}^2.
$$
The fact that this is a well-defined action of $\g$ on $\qq$
translates into the 1-cocycle condition:
$H_{\ga\rho}=(H_{\ga})^{\rho}+H_{\rho}$. Note that $H_{\ga}=0$
precisely when $\ga$ permutes the four lines $x,y,z,x+y+z$.

With this notation, we have $\ga(C_Q)=C_{\ga(Q)}$ and
$$
\op{Isom}(C_Q,C_{Q'})=\{\rho\in\g\tq \rho(Q)=Q'\}, \quad\forall
Q,Q'\in\qq.
$$

Througouht the paper we shall abuse of notation and denote a
quadratic form (\ref{qf}) simply by $Q=(a,b,c,d,e,f)$.

We want to study now the invariants of ordinary quartics. We know
that every ordinary quartic is $\kb$-isomorphic to one of the form
$C_Q$ and we want to study invariants of these curves under the
action of $\g$. In order to obtain a linear action we consider now
models
$$C_Q^{\textrm{pr}} \colon\quad Q^2=g^2 x y z(x+y+z).$$
By considering $(a,b,c,d,e,f,g)$ as homogeneous coordinates, we
can identify the set of these quartics with an open subset of
$\pr{6}$. Consider the following generators of $\g$,
$$
A,\,B\colon\ \pr{2} \longrightarrow \pr{2},\qquad
A(x,y,z)=(x,y,x+z),\ B(x,y,z)=(y,z,x).
$$
Their action on $C_Q^{\textrm{pr}}$ induces a linear
representation on $\pr{6}$ given by the matrices:
$$A = \left(\begin{array}{ccccccc}
1 & 0& 0 & 0 & 0 & 0 & 0 \\
 0 & 1& 0 & 0 & 0 & 0 & 0 \\
1 & 0& 1 & 0 & 0 & 0 & 0 \\
0 & 0& 0 & 1 & 0 & 0 & 0 \\
0 & 0& 0 & 1 & 1 & 0 & 0 \\
1 & 0& 0 & 0 & 0 & 1 & 0 \\
0 & 0& 0 & 1 & 0 & 0 & 1 \\
\end{array}\right), \quad B=\left(\begin{array}{ccccccc}
0 & 1& 0 & 0 & 0 & 0 & 0 \\
 0 & 0& 1 & 0 & 0 & 0 & 0 \\
1 & 0& 0 & 0 & 0 & 0 & 0 \\
0 & 0& 0 & 0 & 1 & 0 & 0 \\
0 & 0& 0 & 0 & 0 & 1 & 0 \\
0 & 0& 0 & 1& 0 & 0 & 0 \\
0 & 0& 0 & 0 & 0 & 0 & 1 \\
\end{array}\right).$$
Changing the model, we may simplify the action. The change of
coordinates $$(a,b,c,d,e,f,g) \mapsto
(a+g,b+g,c+g,d+g,e+g,f+g,g)$$ corresponds to substituting the
models $C_Q$ by models
$$(a x^2 +b y^2 +c z^2+d xy+eyz+f zx)^2=g^2 C_K,$$
with $C_K = x^4+y^4+z^4+(xy)^2+(yz)^2+(zx)^2+x^2 yz+ xy^2z+xy
z^2$. This curve $C_K$ is a twist of the Klein quartic and so it
is invariant by the action of $\g$. Thus, the action of $\g$
restricts to the conic $Q$. Let
 $R=\ft[a,b,c,d,e,f]^{\g}$ be the ring of invariants of conics
under the natural action of $\g$ on $\pr{2}$.

\begin{prop}
The locus of ordinary quartics in $\mm_3$
 is isomorphic to an explicit open subset (given by non-singularity conditions)
 of the affine variety $\op{Spec}(R)$.
\end{prop}

The structure of the ring $R$ is analyzed in more detail in
\cite{mr}.

\subsection{Descent data}\label{ddata}
Let $\sg$ be the Frobenius automorphism of $\kb$ relative to $k$:
$\si a=a^q$, $\forall a\in \kb$.

By Proposition \ref{pr1}, the descents to $k$ of all curves in the
family $\cq$ take into account all $k$-isomorphy classes of
non-singular quartic plane curves defined over $k$.
%In this paragraph we obtain an explicit description of these $k$-isomorphy classes by determining
%all possible descents of the curves $C_Q$.

Let us briefly recall the basic facts of descent theory for the
specific case of curves defined over finite fields. The basic
reference is \cite{w}. Given a curve $\cc$ defined over $\kb$, a
family of descent data for $\cc$ over $k$ is generated by any
isomorphism, $\ga\colon\cc \to \si\cc$, such that
\begin{equation}\label{descent}
\sij{n-1}\ga\circ\cdots\circ\si\ga\circ\ga=1,
\end{equation}for some $n\ge1$, which is called the degree of descent.
We call the pair $(\cc,\ga)$ a {\it descent datum} over $k$. To
such datum we can associate a curve $C$ defined over $k$ and a
$\kb$-isomorphism $\phi\colon C\to \cc$ such that $\ga=\si
\phi\circ\phi^{-1}$; this curve $C$ is unique up to
$k$-isomorphism.
 The degree of descent is the degree of the
minimum field of definition of the isomorphism $\phi$. We denote
by $\dsc(\cc,\ga)$ the class of $k$-isomorphism of $C$. We have:
\begin{equation}\label{eqdd}
\dsc(\cc,\ga)=\dsc(\cc',\ga') \sii \exists \rho\colon \cc\iso \cc'
\mbox{ such that } \ga'=\si\rho\circ\ga\circ\rho^{-1}.
\end{equation}

The descent theory of the family $\cq$ is simplified by the fact
that the isomorphisms involved in the descent data are galois
invariant. In fact, if $(C_Q,\ga)$ is a descent datum over $k$, we
have necessarily $\ga\in\g$, since $\si(C_Q)=C_{\si Q}$.
Therefore, condition (\ref{descent}) amounts in our case to
$\ga^n=1$. Since the elements in $\g$ have order 1,2,3,4 or 7, we
get another proof that these are the possible degrees of descent.

We abuse of language and consider our descent data to be pairs
$(Q,\ga)$ instead of $(C_Q,\ga)$. We denote by $\dd$ the set of
all descent data over $k$ of curves in the family $\cq$:
$$
\dd=\{(Q,\ga)\tq \ga\in\g,\,Q\in\qq,\, \ga(Q)=\si Q\}.
$$If we denote by $\cc_k$ the set of $k$-isomorphy classes of non-singular quartic plane
curves whose Jacobian is ordinary, we get an onto map $\dsc\colon
\dd \lra \cc_k$. By (\ref{eqdd}):
\begin{equation}\label{prop2}
\dsc(Q,\ga)=\dsc(Q',\ga') \sii \exists \rho\in\g \mbox{ such that
}\rho(Q)=Q'\mbox{ and } \ga'=\rho\ga\rho^{-1}.
\end{equation}

%\begin{proof}
%Let $\phi$, resp. ${\phi'}$, be $K$-isomorphisms $\phi\colon C\to C_Q$,  resp. ${\phi'}\colon C' \to
%C_{Q'}$, such that $\ga=\si \phi\circ\phi^{-1}$, resp. $\ga'=\si {\phi'}\circ{\phi'}^{-1}$.
%If $\vf\colon C \lra C'$ is a $k$-isomorphism, then $\rho:=\phi'\circ \vf \circ \phi^{-1}$ is an
%isomorphism from $C_Q$ to $C_{Q'}$. Hence, $\rho$ belongs to $\g$ and one checks easily that
%$\ga'=\rho \ga \rho^{-1}$.
%Conversely, assume that $\ga'=\rho \ga \rho^{-1}$ for some $\rho\in\g$ such that $\rho(Q)=Q'$; then,
%$\vf:={\phi'}^{-1} \circ \rho \circ \phi$ is an isomorphism between $C$ and $C'$ and
%$$\si \vf=\si{\phi'}^{-1} \circ \rho\circ \si\phi=({\phi'}^{-1} \circ {\ga'}^{-1}) \circ \rho \circ
%(\ga \circ \phi)=\vf.$$ Thus, $C$ and $C'$ are isomorphic over $k$.
%\end{proof}

\begin{lem}\label{bg}
Let $C$, $C'$ be ordinary curves defined over $k$, which have been
obtained by descent from respective pairs $(Q,\ga)$, $(Q',\ga)$,
with the same $\ga\in\g$. Then, the sets $\bb$, $\bb'$ of
bitangents of $C$ and $C'$ are $\pg{k}$-equivalent.
\end{lem}

\begin{proof}
Let $\phi$, resp. $\phi'$, be $\kb$-isomorphisms $\phi\colon C\to
C_Q$,  resp. $\phi'\colon C' \to C_{Q'}$, such that
$\si\phi\circ\phi^{-1}=\ga=\si \phi'\circ{\phi'}^{-1}$. We have
$\phi(\bb)= \bb_0=\phi'(\bb')$. Hence,
$\phi^{-1}\circ\phi'(\bb')=\bb$, and clearly $\phi^{-1}\circ\phi'$
is defined over $k$.
\end{proof}\e

\ni{\bf Notation. }We denote by $\gga:= \{\rho\in\g\tq
\rho\ga=\ga\rho\}$  the centralizer of $\ga$ in $\g$. Moreover,
for any subset $\bb$ of $\pr{2}$ and field $K\subseteq \kb$ we
denote
$$\aut_K(\bb):=\{\rho\in\pg{K}\tq \rho(\bb)=\bb\}.$$For $K=\kb$ we write simply $\aut(\bb)$.

\begin{lem}\label{aut}
Let $C$ be an ordinary curve defined over $k$ which has been
obtained by descent from $(Q,\ga)$ and let $\phi\colon C\iso C_Q$
be an isomorphism such that $\si\phi\circ\phi^{-1}=\ga$. Then, for
the set of bitangents $\bb$ of $C$ we have:
$$\aut_k(\bb)=\phi^{-1}\circ\gga\circ\phi.$$

Assume moreover that $C'$ is another ordinary curve such that
$\phi(C')\in\cq$, say $\phi(C')=C_{Q'}$. Then,
$$
\op{Isom}_k(C,C')=\phi^{-1}\circ\{\rho\in\gga\tq
\rho(Q)=Q'\}\circ\phi.
$$In particular, $ \aut_k(C)=\phi^{-1}\circ \gga(Q)\circ\phi$, where $\gga(Q)$ denotes the
stabilizer of $Q$ under the action of $\gga$.
\end{lem}

\renewcommand\arraystretch{1.4}
\begin{proof}
We have clearly:
$$
\begin{array}{l}
\aut(\bb)=\phi^{-1}\circ \aut(\bb_0)\circ \phi=\phi^{-1}\circ \g\circ \phi,\\
\op{Isom}(C,C')=\phi^{-1}\circ \op{Isom}(C_Q,C_{Q'})\circ
\phi=\phi^{-1}\circ \{\rho\in\g\tq \rho(Q)=Q'\}\circ \phi.
\end{array}
$$Finally, one checks immediately that for any $\rho\in\g$ the automorphism $\phi^{-1}\circ\rho
\circ\phi$ is defined over $k$ if and only if $\rho$ commutes with
$\ga$.
\end{proof}
\renewcommand\arraystretch{1.}

For any fixed $\ga\in\g$, let us denote by $\ddg$ the set of all
pairs $(Q,\ga)$ belonging to $\dd$. By Lemma \ref{bg}, the sets of
bitangents of all curves obtained by descent of all elements in
$\ddg$ lie in the same orbit under the action of $\pg{k}$; we
denote by $\bb_{\ga}$ this common orbit.

\begin{prop}\label{prop3}
For any $\ga,\ga'\in \g$, the following conditions are equivalent:
\begin{enumerate}
\item There exists $\rho \in \g$ such that $\ga' = \rho \ga \rho^{-1}$.
\item $\dsc(\ddg)=\dsc(\dd_{\ga'})$.
\item $\dsc(\ddg) \cap \dsc(\dd_{\ga'}) \ne \emptyset$.
\item $\bb_{\ga}=\bb_{\ga'}$.
\end{enumerate}
\end{prop}

\begin{proof}
$1$ implies $2$: Since $\rho$ is galois invariant, for any descent
datum $(Q,\ga)\in\ddg$, we see that $(\rho(Q),\ga')$ is a descent
datum in $\dd_{\ga'}$. By (\ref{prop2}),
$\dsc((Q,\ga))=\dsc((\rho(Q),\ga'))$; hence $\ddg\subseteq
\dd_{\ga'}$ and by symmetry these two sets coincide.

$2$ implies $3$ and $3$ implies $4$ are trivial.

$4$ implies $1$: Let $\phi\colon C\to C_Q$, resp. $\phi'\colon C'
\to C_{Q'}$, be isomorphisms such that $\ga=\si
\phi\circ\phi^{-1}$, resp. $\ga'=\si \phi'\circ{\phi'}^{-1}$.

Let $\bb$, $\bb'$ be the respective sets of bitangents of $C$,
$C'$ and assume that there is $\eta\in\pg{k}$ such that
$\eta(\bb)=\bb'$. Then, $\rho:=\phi'\circ\eta\circ\phi^{-1}$
belongs to $\g$, since it leaves $\bb_0$ invariant. Moreover,
$\rho\ga\rho^{-1}=\si\rho\ga\rho^{-1}=\ga'$.
\end{proof}

The centralizer  $\gga$ of any $\ga\in\g$ operates on $\ddg$ by:
$\rho(Q,\ga)=(\rho(Q),\ga)$. By (\ref{prop2}), two descent data in
$\ddg$ are equivalent iff they are in the same orbit under this
action. Therefore, we have proved the following:

\begin{teor}\label{bij}
For any system of representatives $\cc(\g)$ of conjugacy classes
of $\g$, the mapping $\dsc$ determines a bijection
$$\qquad\qquad\dsc \colon \amalg_{\ga\in\cc(\g)} \; \gga\backslash\ddg  \lra \cc_k.\qquad\qquad\qed $$
\end{teor}

All elements of order 1,2,3,4 of $\g$ lie respectively in one
single conjugacy class, whereas the elements of order 7 split into
two conjugacy classes, according to their trace being 0 or 1
\cite[Section 7.4]{h}. We fix once and for all,
\begin{equation}\label{cc}
\cc(\g)=\{1,\,\ga_2,\,\ga_3,\,\ga_4,\,\ga_{7,0},\,\ga_{7,1}\},
\end{equation} where:
$$\ga_2=\begin{pmatrix} 0&1&0\\1&0&0\\0&0&1\end{pmatrix},\
\ga_3=\begin{pmatrix} 0&1&0\\0&0&1\\1&0&0\end{pmatrix},\
\ga_4=\begin{pmatrix} 0&1&0\\0&0&1\\1&1&1\end{pmatrix},
$$$$
\ga_{7,0}=\begin{pmatrix}0&1&0\\0&0&1\\1&1&0\end{pmatrix},\
\ga_{7,1}=\begin{pmatrix} 0&1&0\\0&0&1\\1&0&1\end{pmatrix}.
$$

By Theorem \ref{bij}, we get a splitting of $\cc_k$ into the disjoint union of six subfamilies:
$$
\cc_k=\cc_1\amalg\cc_2\amalg\cc_3\amalg\cc_4\amalg\cc_{7,0}\amalg\cc_{7,1},
$$ where the subindex indicates the degree of the descent, which by Lemma \ref{degdd} is the degree
of the minimum field of definition of all bitangents. By
Proposition \ref{prop3} these families are in one-to-one
correspondence with the orbits of Fano planes defined over $k$ (as
a set) under the action of $\pg{k}$.

In the next two paragraphs we compute the cardinalities of each of
these families and we exhibit rational normal models for them.

\subsection{Number of ordinary curves}\label{nbr}
We denote by $k_2$, $k_3$, $k_4$, $k_7$ the respective extension
of $k$ in $\kb$, of degree 2,3,4,7. The action of Frobenius will
be denoted sometimes by $(\ )'$: $\si a=a'$, $\sij{2}a=a''$, etc.

After Theorem \ref{bij}, in order to compute $\cc_k$ we need only
to compute $|\gga\backslash\ddg|$ for $\ga$ running on the system
of representatives of conjugacy classes of $\g$ specified in
(\ref{cc}).

In general, when we have a finite group $G$ acting on a finite set
$X$, the number of orbits can be
 counted with the well-known formula:
\begin{equation}\label{comb}
\vert G \backslash X\vert = \frac{1}{\vert G \vert} \sum_{g\in G}
\vert X(g)\vert=\sum_{g\in\cc(G)}\frac{|X(g)|}{|G_g|},
\end{equation}
where $X(g):=\{x\in X \tq g(x)=x\}$ is the set of fixed points of
$g$, $\cc(G)$ is a system of representatives of conjugacy classes
of $G$ and $G_g$ is the centralizer of $g$.

\begin{teor} There are $q^6-q^5+q^4-3q^3+5q^2-6q+7$
$\fq$-isomorphism classes of ordinary non-hyperelliptic curves of
genus three defined over $\fq$. More precisely, according to the
minimum field of definition of its bitangents, these curves are
distributed in the following way:\be

\renewcommand\arraystretch{1.4}
\centerline{\begin{tabular}{|c|c|} \hline
$\cc_1$&$\frac1{168}\left(q^6-7q^5+42q^4-140q^3+343q^2-462q+328\right)$\\
\hline $\cc_2$&$\frac
18\left(q^6-3q^5+6q^4-12q^3+15q^2-6q\right)$\\\hline
$\cc_3$&$\frac
13\left(q^6-q^5-2q^3+4q^2-6q+7\right)$\\
\hline $\cc_4$&$\frac 14\left(q^6-q^5-q^2-2q+4\right)$\\
\hline $\cc_{7,0}+\cc_{7,1}$&$\frac 17(q^6+6)+\frac 17(q^6+6)$\\
\hline
\end{tabular}}\be
\end{teor}
\renewcommand\arraystretch{1.}

\begin{proof}
For commodity, we drop the $\ga$ from the pair $(Q,\ga)$. More
precisely, we identify $\ddg=\{Q\in\qq\tq \ga(Q)=\si Q\}$, for any
$\ga\in\g$.\e

\ni{\bf Case 1. }For $\ga=1$ we have $\gga=\g$ and
$\ddg=\{Q\in\qq\tq Q=\si Q\}=:\qq_k$ coincides with the set of
quadratic forms in $\qq$ with coefficients in $k$.

The cardinality of the centralizers of $\g$ is well-known
\cite{h}: $|\g_{\rho}|=168,\,8,\,3,\,4,\,7$, according to the
order of $\rho$ being respectively 1,2,3,4,7. Thus, in order to
apply (\ref{comb}) we need only to compute $|\qq_k(\rho)|$ for all
$\rho\in\cc(\g)$.

\begin{lem}\label{abne1}
$|\{(a,b)\in k^*\times k^*\tq a+b\ne 1\}|=q^2-3q+3$.
\end{lem}

\begin{proof}
For $a=1$, resp. $a\ne 1$, we have $q-1$, resp. $q-2$
possibilities for $b$.
\end{proof}

Let $\uu=(k^*)^3\times \,k^3$. In order to compute $|\qq_k|$ we
need to count how many elements $(a,b,c,d,e,f)\in\uu$ fail to
satisfy at least one of the following linear relations:
$$
a+b+d=0,\ b+c+e=0,\ a+c+f=0,\ a+b+c+d+e+f=1.
$$

For $i=1,2,3$, there are $(q-1)^3q^{3-i}$ elements in $\uu$
satisfying exactly $i$ of these relations, since we can fix any
value of $(a,b,c)\in (k^*)^3$ and $i$ of the values of $d,e,f$ are
determined by the linear relations. By Lemma \ref{abne1}, there
are $q^2-3q+3$ elements in $\uu$ satisfying all four linear
equations simultaneously, since we can choose $(a,b)\in (k^*)^2$
with $a+b\ne1$ arbitrarily and then determine $c,d,e,f$ by the
relations.

By the inclusion-exclusion principle, we have
\begin{multline*}
|\qq_k|=|\uu|-4(q-1)^3q^2+6(q-1)^3q-4(q-1)^3+(q^2-3q+3)=\\
=q^6-7q^5+21q^4-35q^3+35q^2-21q+7.
\end{multline*}
On the other hand,
$$
\qq_k(\ga_2)=\{(a,a,c,d,e,e)\in k^6\tq acd\ne0,\, a+c+e\ne0,\,
c+d\ne 1\}.
$$
By Lemma \ref{abne1}$, |\qq_k(\ga_2)|=(q^2-3q+3)(q-1)^2$, since we
can choose $(c,d)\in (k^*)^2$ with $c+d\ne 1$ and then we have
$q-1$ possibilities, both for $a$ and for $e$. Moreover,
\renewcommand\arraystretch{1.4}
$$
\begin{array}{lll}
\qq_k(\ga_3)=\{(a,a,a,d,d,d)\in (k^*)^6\tq a+d\ne 1\}&\imp&
|\qq_k(\ga_3)|=q^2-3q+3,\\
\qq_k(\ga_4)=\{(a,b,a,b,b,b)\in (k^*)^6\}& \imp& |\qq_k(\ga_4)|=(q-1)^2,\\
\qq_k(\ga_{7,0})=\qq_k(\ga_{7,1})=\{(1,1,1,1,1,1)\}& \imp&
|\qq_k(\ga_{7,0})|=|\qq_k(\ga_{7,1})|=1.
\end{array}
$$
\renewcommand\arraystretch{1.}
By (\ref{comb}), we get:
$$
|\g\backslash\qq_k|=\frac
1{168}\,\left(q^6-7q^5+42q^4-140q^3+343q^2-462q+328\right).
$$\e

\ni{\bf Case 2. } For $\ga=\ga_2$ we have
$\gga=\{1,\,\ga,\,\tau,\,\ga\tau,\,\rho,\,\rho^3,\,\rho\,\tau,\,\tau\rho\}\simeq
D_8$, where
\begin{equation}\label{d8}
\tau=\begin{pmatrix} 1&0&0\\0&1&0\\1&1&1\end{pmatrix},\
\rho=\begin{pmatrix} 1&0&1\\0&1&1\\1&1&1\end{pmatrix}.
\end{equation}
The center of $\gga$ is $\{1,\ga\}$ and $\tau^2=1$, $\rho^2=\ga$.
The conjugacy classes of $\gga$ are:
$$
\gga=\{1\}\amalg\{\ga\}\amalg\{\tau,\,\ga\tau\}\amalg\{\rho,\rho^3\}\amalg
\{\rho\,\tau,\,\tau\rho\}.
$$
We have now,
\begin{multline*}
\ddg=\{(a,a',c,d,e',e)\in (k_2)^6\tq ac\ne0,\, c,d\in k,\\
d\ne\tr(a),\,a+c+e\ne0,\,c+d+\tr(a)+\tr(e)\ne1\}.
\end{multline*}

Let $\uu=\{(a,a',c,d,e',e)\in (k_2)^6\tq ac\ne0,\, c,d\in k\}$,
with $|\uu|=(q^2-1)(q-1)q^3$.

We want to count how many elements in $\uu$ fail to satisfy at
least one of the following linear relations:
$$
d=\tr(a),\ a+c+e=0,\ c+d+\tr(a)+\tr(e)=1.
$$

There are $(q^2-1)(q-1)q^2$ elements in $\uu$ satisfying
$d=\tr(a)$ or $c+d+\tr(a)+\tr(e)=1$, whereas only $(q^2-1)(q-1)q$
elements satisfy $e=a+c$.  There are $(q^2-1)(q-1)$ elements
satisfying $e=a+c$ and one of the other two relations, whereas
$(q^2-1)(q-1)q$ elements satisfy $d=\tr(a)$ and $c+\tr(e)=1$.
Finally, there are $q^2-q-1$ elements satisfying simultaneously
the three equations (necessarily $\tr(a)\ne 1$ and the values of
$c,d,e$ are determined by the choice of $a$).

By the inclusion-exclusion principle, we have
\begin{multline*}
|\ddg|=|\uu|-(q^2-1)(q-1)(2q^2+q)+(q^2-1)(q-1)(q+2)-(q^2-q-1)=\\
=q^6-3q^5+q^4+5q^3-5q^2-q+3.
\end{multline*}
Now, $\ddg(\ga)=\qq_k(\ga_2)$ has $(q^2-3q+3)(q-1)^2$ elements, as
we saw in Case 1, whereas,
\renewcommand\arraystretch{1.4}
$$\begin{array}{l}
\ddg(\rho)=\{(a,a,c,1,c,c)\in(k^*)^6\}\ \imp |\ddg(\rho)|=(q-1)^2.\\
\ddg(\tau)=\{(a,a',c,d,c,c)\in (k_2)^6\tq ac\ne0,\,c,d\in
k,\,d\ne\tr(a),\,
c+d+\tr(a)\ne1\},\\
\ddg(\rho\,\tau)=\{(a,a,c,e+e'+1,e',e)\in (k^*)^4\times (k_2)^2\tq
\tr(e)\ne1,\, a+c+e\ne0\}.
\end{array}
$$\renewcommand\arraystretch{1.}
We have $|\ddg(\tau)|=(q^2-1)(q^2-3q+3)$, since for any fixed
$a\in (k_2)^*$ there are $q^2-3q+3$ possibilities for the pair
$(c,d+\tr(a))$ (hence for the pair $(c,d)$) by Lemma \ref{abne1}.
Moreover, $|\ddg(\rho\,\tau)|=(q^2-q-1)(q-1)^2$ since for any
fixed pair $(a,c)\in (k^*)^2$ there are $q^2-q-1$ possible values
for $e$. We get, finally, by (\ref{comb}):
$$
|\gga\backslash\ddg|=\frac
18\,\left(q^6-3q^5+6q^4-12q^3+15q^2-6q\right).
$$\e

\ni{\bf Case 3. } For $\ga=\ga_3$ we have
$\gga=\{1,\,\ga,\,\ga^2\}\simeq C_3$ and
$$
\ddg=\{(a,a',a'',d,d',d'')\in (k_3)^6\tq a\ne 0,\, a+a'+d\ne
0,\,\tr(a)+\tr(d)\ne 1\}.
$$The set $\uu:=k_3^*\times k_3$ has $(q^3-1)q^3$ elements. There
are $q^3-1$ elements $(a,d)\in\uu$ satisfying $d=a+a'$,
$(q^3-1)q^2$ elements satisfying $\tr(d)=1+\tr(a)$ and $q^2$
elements satisfying both conditions (necessarily $\tr(a)=1$).
Hence,
$$
|\ddg|=|\uu|- (q^3-1)(q^2+1)+q^2=q^6-q^5-2q^3+2q^2+1.
$$
For $i=1,2$, $\ddg(\ga^i)=\qq_k(\ga_3)$ has $q^2-3q+3$ elements,
as we saw in Case 1. We get:
$$
|\gga\backslash\ddg|=\frac
13\,\left(q^6-q^5-2q^3+4q^2-6q+7\right).
$$\e

\ni{\bf Case 4. } For $\ga=\ga_4$ we have
$\gga=\{1,\,\ga,\,\ga^2,\,\ga^3\}\simeq C_4$ and
$$
\ddg=\{(c',b,c,b+c'+c'',b+c+c''',b'+c+c')\tq c\in k_4^*,\, b\in
k_2^*,\,\tr_{k_4/k}(c)+\tr_{k_2/k}(b)\ne 1\}.
$$ There are $q^3$ elements $c\in k_4^*$ with
$\tr_{k_4/k}(c)=1$; for each of them, we can choose $q^2-q$
elements $b\in k_2^*$ with $\tr_{k_2/k}(b)\ne 0$. There are
$q^4-q^3-1$ elements $c\in k_4^*$ with $\tr(c)_{k_4/k}\ne1$; for
each of them, we can choose $q^2-q-1$ elements $b\in k_2^*$ with
$\tr_{k_2/k}(b)\ne 1+\tr_{k_4/k}(c)$. Altogether, we have
$$
|\ddg|=q^3(q^2-q)+(q^4-q^3-1)(q^2-q-1)=q^6-q^5-q^4+q^3-q^2+q+1.
$$

For $i=1,3$, $\ddg(\ga^i)=\qq_k(\ga_4)$ has $(q-1)^2$ elements, as
we saw in Case 1, whereas,
\begin{multline*}
\ddg(\ga^2)=\{(c',b,c,b+c+c',b+c+c',b'+c+c')\tq
b,c\in k_2^*,\, \tr_{k_2/k}(b)\ne 1\}\ \imp\\
|\ddg(\ga^2)|=(q^2-1)(q^2-q-1)=q^4-q^3-2q^2+q+1.
\end{multline*}
We get:
$$
|\gga\backslash\ddg|=\frac 14\,\left(q^6-q^5-q^2-2q+4\right).
$$\e

\ni{\bf Cases (7,0), (7,1). } For $\ga=\ga_{7,0}$ we have
$\gga=\gen{\ga}\simeq C_7$ and
$$
\ddg=\{(b'',b,b',b+b''+b''',b+b'+b^{v},b'+b''+b^{vi})\tq b\in
k_7^*,\, \tr(b)=1\}.
$$ Now, $|\ddg|=q^6$, $\ddg(\ga^i)=\{(1,1,1,1,1,1)\}$
for $i=1,\dots,6$, and
$$
|\gga\backslash\ddg|=\frac 17\,\left(q^6+6\right).
$$
For $\ga=\ga_{7,1}$ we obtain a completely analogous result.
\end{proof}

\subsection{Rational normal models}\label{ratmod}
In paragraph \ref{ddata} we obtained a partition:
$$
\cc_k=\cc_1\amalg\cc_2\amalg\cc_3\amalg\cc_4\amalg\cc_{7,0}\amalg\cc_{7,1},
$$
of the $k$-isomorphy classes of ordinary plane quartics, according
to the minimum field of definition of their set of bitangents. In
this section we exhibit rational models for each of the six
families.

By Theorem \ref{bij}, each family is identified to
$\g_{\ga}\backslash\ddg$ for some $\ga\in\cc(\g)$; moreover, by
Proposition \ref{prop3} it corresponds to an orbit of Fano planes
defined over $k$ (as a set) under the action of $\pg{k}$. We shall
choose a Fano plane $\bb$ in each orbit and consider the family
$\oo$ of normal models:
$$N_Q\colon\qquad Q^2=\ell_1\ell_2\ell_3(\ell_1+\ell_2+\ell_3),$$ where
$\ell_1,\ell_2,\ell_3$ are three fixed non-concurrent lines of
$\bb$ and $Q$ is a quadratic form such that $N_Q$ is defined over
$k$ and non-singular. The nonsingularity condition amounts to
$N_Q$ not passing through any of the seven points of the Fano
plane $\bb$.

In all cases, the automorphism $\phi$ of $\pr{2}$ given by:
$$
\phi(x,y,z)=(\ell_1(x,y,z),\ell_2(x,y,z),\ell_3(x,y,z)),
$$satisfies $\si\phi\circ \phi^{-1}=\ga$ and the mapping
$$
\oo\lra \ddg,\qquad N_Q\mapsto (\phi(N_Q),\ga)
$$establishes a bijection between $\oo$ and $\ddg$. Thus, the models $N_Q$ faithfully represent
$\dsc(\ddg)$. Moreover, by Lemma \ref{aut} the group
$\aut_k(\bb)=\phi^{-1}\circ\gga\circ\phi$ acts fully by
$k$-isomorphisms on the family $\oo$.

We shall choose the first three lines of $\bb$ in such a way that
$\phi$ will be represented by a symmetric matrix. Then,
$$
\phi^{-1}\circ \ga\circ\phi=\phi^{-1}\circ \si\phi=\t\ga,
$$and $\aut_k(\bb)=\g_{\t\ga}$.
As in section \ref{action}, the action of $\g_{\t\ga}$ on the
family of the curves $N_Q$ can be interpreted in terms of an
action on the corresponding set of quadratic forms $Q$. This
action has again the shape: $\rho(Q)=Q^{\rho^{-1}}+H_{\rho^{-1}}$,
where the term $H_{\rho^{-1}}$ depends on how $\rho$ permutes the
seven lines of $\bb$. More precisely, $H_{\rho}$ is defined by:
$$
\ell_1^{\rho}\ell_2^{\rho}\ell_3^{\rho}(\ell_1+\ell_2+\ell_3)^{\rho}=
\ell_1\ell_2\ell_3(\ell_1+\ell_2+\ell_3)+H_{\rho}^2.
$$The explicit description of this action allows us to determine the $k$-isomorphy classes
contained in $\oo$ and to compute the $k$-automorphism group of
each curve $N_Q$.

Any nonsingular quartic plane curve $C$ is $k$-isomorphic to one
of our models. By computing its set of bitangents we know to what
family it corresponds; any $k$-automorphism of $\pr{2}$ taking the
set of bitangens of $C$ to $\bb$ will furnish a $k$-isomorphism of
the curve $C$ with one of the normal models.

\subsubsection*{Normal models $\oo_1$}
We choose $\bb=\bb_0=\{x,\, y,\,z,\,x+y+z,\,x+y,\,y+z,\,x+z\}$,
with $\aut_k(\bb)=\g$. The family $\oo_1$ gathers the curves
$N_Q:=C_Q$ with $Q$ defined over $k$:
$$N_Q=C_Q\colon\quad Q^2=xyz(x+y+z), \qquad Q\in \qq_k.$$

The action of $\g$ on $\qq_k$ is the ordinary action that we
introduced in section \ref{action}. Two curves $C_Q$, $C_{Q'}$ in
this family are $k$-isomorphic iff they are $\kb$-isomorphic iff
$Q'=\rho(Q)$ for some $\rho\in\g$. The different possibilities for
the group $\aut_k(C_Q)=\aut(C_Q)$ in terms of $Q$ can be found in
 \cite{wa}.\renewcommand\arraystretch{1.2}

%\footnote{We could also have chosen the product of any other
%family of three linearly independent forms in $\bb$ and their sum
%and adapt the conditions on $Q$ to ensure non-singularity.}
\subsubsection*{Normal models $\oo_2$}
We fix as a generator of $k_2/k$ an element $u\in k_2\setminus k$,
 with equation $u^2+u=r$, for certain $r\not\in\as(k):=\{x+x^2\tq x\in k\}$. We
choose  $\ell(x,y,z)=ux+u'y$ and
$$
\begin{array}{l}
\bb=\{\ell,\,\ell',\,z,\,x+y+z,\,x+y,\,\ell'+z,\,\ell+z\},\\
%c\ne0,\,a+b+d\ne0,\,a+b+c+d+e+f\ne1,\,(a+b,a+dr)\ne(0,0),\,
%(a+b+e+f,a+c+f+r(a+b+d))\ne (0,0)\}.
%\begin{multline*}
 \qq_2:=\{Q(x,y,z)=ax^2+by^2+cz^2+dxy+eyz+fxz\tq a,b,c,d,e,f\in
k,\qquad\quad\\\quad\
c\ne0,\,a+b+d\ne0,\,a+b+c+d+e+f\ne1,\,Q(u,u',0)\ne0,\,Q(u,u',1)\ne0\}.
\end{array}$$
\renewcommand\arraystretch{1.}
%\end{multline*}
We consider the family $\oo_2$ of normal models:
$$
N_Q\colon\quad Q^2=\ell\ell'z(x+y+z)=(rx^2+ry^2+xy)z(x+y+z),\qquad
Q\in\qq_2.
$$
Since $\ga:=\ga_2$ is symmetric, we have
$$\aut_k(\bb)=\g_{\t\ga}=\gga=\{1,\,\ga,\,\tau,\,\ga\tau,\,\rho,\,\rho^3,\,\rho\,\tau,\,\tau\rho\}
\simeq D_8,$$ where  $\tau,\rho$ are given in (\ref{d8}). The
orbit of $Q=(a,b,c,d,e,f)$ under the action of $\gga$ is
$$
\begin{array}{ccccccc}
a&b&c&d&e&f&\qquad (1)\\
b&a&c&d&f&e&\qquad (\ga)\\
a+c+f&b+c+e&c&d+e+f&e&f&\qquad (\tau)\\
b+c+e&a+c+f&c&d+e+f&f&e&\qquad (\ga\tau)\\
b+c+e&a+c+f&\Sigma+1&d+e+f&d+e+1&d+f+1&\qquad (\rho)\\
a+c+f&b+c+e&\Sigma+1&d+e+f&d+f+1&d+e+1&\qquad (\rho^{-1})\\
b&a&\Sigma+1&d&d+f+1&d+e+1&\qquad (\rho\,\tau)\\
a&b&\Sigma+1&d&d+e+1&d+f+1&\qquad (\tau\rho),
\end{array}
$$
where $\Sigma=a+b+c+d+e+f$. Hence,
$$
\aut_k(N_Q)=\left\{
\begin{array}{ll}
\gga\simeq D_8,&\mbox{ if }a=b,\,d=1,\,c=e=f,\\
\{1,\ga,\tau\rho,\rho\,\tau\}\simeq C_2\times C_2,&\mbox{ if }a=b,\,d=1,\,e=f\ne c,\\
\{1,\ga,\tau,\ga\tau\}\simeq C_2\times C_2,&\mbox{ if }a=b,\,d\ne1,\,c=e=f,\\
\gen{\tau\rho}\simeq C_2,&\mbox{ if }a=b+e+f,\,d=1,\,e\ne f,\\
\gen{\rho\,\tau}\simeq C_2,&\mbox{ if }a=b,\,d+e+f=1,\,e\ne f,\\
\gen{\tau}\simeq C_2,&\mbox{ if }c=e=f,\,a\ne b\\
\gen{\ga\tau}\simeq C_2,&\mbox{ if }e=f=a+b+c,\,a\ne b\\
\gen{\ga}\simeq C_2,&\mbox{ if }a=b,\,e=f\ne c,\,d\ne1\\
\{1\},&\mbox{ otherwise}.
\end{array}
\right.
$$\renewcommand\arraystretch{1.2}

\subsubsection*{Normal models $\oo_3$}
We fix as a generator of $k_3/k$ an element $v\in k_3\setminus k$,
with equation $v^3+v^2=s$. By \cite[Lemma 7]{cnp}, we have
$v'=v^2(tv+1)^{-1}$, $v''=v^2((t+1)v+1)^{-1}$, where $t\in k$
satisfies $t^2+t+1=s^{-1}$. We take $\ell(x,y,z)=vx+v'y+v''z$ and
$$
\begin{array}{l}
\bb=\{\ell,\,\ell',\,\ell'',\,x+y+z,\,\ell+\ell',\,\ell'+\ell'',\,\ell+\ell''\},\\
%\begin{multline*}
\qq_3:=\{Q=(a,b,c,d,e,f)\in k^6\tq \qquad\qquad \qquad\qquad\qquad
\qquad\qquad \qquad\qquad
\qquad\quad \\
 \qquad\qquad \qquad\qquad Q(v,v',v'')\ne0,\, Q(v+1,v'+1,v''+1)\ne0,\,Q(1,1,1)\ne1\}.
%\end{multline*}
\end{array}$$\renewcommand\arraystretch{1.}
Each $Q\in\qq_3$ provides a normal model:
\begin{multline*}
N_Q\colon\quad Q^2=\ell\ell'\ell''(x+y+z)=(x+y+z)\cdot\\
\cdot\left(s(x^3+y^3+z^3)+xyz+st(xy^2+x^2z+yz^2)+s(t+1)
(xz^2+x^2y+y^2z)\right).
\end{multline*}

Let us denote $\ga:=\ga_3$. Since $\t\ga=\ga^{-1}$, we have
$\aut_k(\bb)=\gga=\{1,\ga,\ga^2\}$. Since $\ga$ permutes the four
bitangents $\ell,\ell',\ell'',x+y+z$, the action of $\gga$ on the
normal models is given by $\rho(Q)=Q^{\rho^{-1}}$. Thus, the orbit
of $Q=(a,b,c,d,e,f)$ is the cyclic orbit generated by
$\ga(Q)=(b,c,a,e,f,d)$ and all $k$-automorphism groups are trivial
except for:
$$
\aut_k(N_Q)=\gga\simeq C_3,\quad\mbox{ if }a=b=c,\,d=e=f.
$$

\subsubsection*{Normal models $\oo_4$}
We fix as a generator of $k_4/k$ an element $w\in k_4\setminus
k_2$, with equation $w^4+(t+t^2)w^2+t^2w=1$, where $t$ is any
element in $k$ such that $t^{-1}\not\in\as(k)$ \cite[Proposition
6]{ns}. For instance, we can choose $t=1$ if $q$ is not a square.
It is easy to check that $w+w''=t$ and $\al:=w+w'$ satisfies the
equation $\al^2+t\al=t$. We take,
$$
\bb=\{\ell,\,\ell',\,\ell'',\,\ell+\ell'+\ell'',\,\ell+\ell',\,\ell'+\ell'',\,\ell+\ell''\},
$$where $\ell(x,y,z)=wx+w'y+w''z$. Note that $\ell+\ell'+\ell''=\ell'''$ since $w$ has null trace.
The linear forms $\ell+\ell'$, $\ell'+\ell''$ are defined over
$k_2$ and conjugate, whereas $\ell+\ell''=t(x+y+z)$ is defined
over $k$.

The seven points of the Fano plane $\bb$ are $(w+t+1,\al,w+\al+1)$
and its three conjugates, $(\al,t,\al')$ and its conjugate, and
$(1,0,1)$. Hence, we consider
\begin{multline*}
\qq_4:=\{Q=(a,b,c,d,e,f)\in k^6\tq \\
Q(w+t+1,\al,w+\al+1)\ne0,\, Q(\al,t,\al')\ne0,\, Q(1,0,1)\ne
t^2\}.
\end{multline*}
We get as normal models:
$$
N_Q\colon\quad Q^2=\ell\ell'\ell''\ell''',\qquad Q\in\qq_4.
$$

For $\ga:=\ga_4$ we have
$\aut_k(\bb)=\g_{\t\ga}=\{1,\t\ga,\t\ga^2,\t\ga^3\}$. Since
$\t\ga$ permutes the four bitangents $\ell,\ell',\ell'',\ell'''$,
the orbit of $Q=(a,b,c,d,e,f)$ is the cyclic orbit generated by
$\t\ga(Q)= Q^{\t\ga^{-1}}=(a+b+c+d+e+f,a,b,d+f,d,d+e)$ and
$$
\aut_k(N_Q)=\left\{
\begin{array}{ll}
\gen{\t\ga}\simeq C_4,&\mbox{ if }a=b=c,\,f=0,\,d=e,\\
\gen{\t\ga^2}\simeq C_2,&\mbox{ if }a=c,\,d+e+f=0,\,(f\ne0 \mbox{ or }b\ne c),\\
\{1\},&\mbox{ otherwise}.
\end{array}
\right.
$$

\subsubsection*{Normal models $\oo_{7,0}$ and $\oo_{7,1}$}
We are interested in elements $\z\in k_7\setminus k$ such that
$$
\{0\}\cup\{\z,\,\z',\,\z'',\,\z''',\,\z^{iv},\,\z^v,\,\z^{vi}\}
$$is an additive subgroup of $k_7$. If $f(x)\in k[x]$ is the monic minimal polynomial of $\z$
over $k$, this condition is equivalent to $xf(x)$ being an
additive polynomial; hence, it is equivalent to
$f(x)=x^7+ax^3+bx+c$, for some $a,b,c\in k$. Let us denote by
$S\subseteq k[x]$ the set of all septic
 irreducible polynomials of this type. Among these elements $\z$ we can distinguish two cases:
$$
\z+\z'''=\z',\qquad \z+\z'''=\z'',
$$which we denote respectively as ``case 0'' and ``case 1''.
It is easy to check that there are no other possibilities. In
terms of the minimal polynomial we distinguish the two cases
according to $f(x)$ dividing $x^{q^3}+x^q+x$ (case 0) or
$x^{q^3}+x^{q^2}+x$ (case 1). Thus, our set $S$ splits as the
disjoint union $S=S_0\cup S_1$ of two subsets gathering the
irreducible polynomials of each type. We want to ensure that these
subsets $S_0$, $S_1$ are non-empty. Actually, we have:

\begin{lem}\label{nset}
$(x^{q^3}+x^q+x)(x^{q^3}+x^{q^2}+x)=x^2\prod_{f(x)\in S}f(x)$.

In particular, $|S_0|=|S_1|=(q^3-1)/7$.
\end{lem}

\begin{proof}
%A polynomial $x^7+ax^3+bx+c$ is separable iff $c\ne0$. Hence, we have $q^3-q^2$ separable polynomials of
%this type with coefficients in $k$. That is, there are $q^3-q^2$ subgroups $H$ of $\kb$ with
%$\dm_{\ft}(H)=3$ and $H$ defined over $k$. It is easy to count how many of these subgroups admit orbits
%under the galois action that are contained in fields of degree less than seven. If one of these subgroups
%is contained in the extension of degree 5,6,10 or 12 of $k$, it is necessarily contained in a smaller
%extension (exercise), whereas the number of subgroups minimally contained in $k,k_2,k_3,k_4$ is:\be
%\renewcommand\arraystretch{1.4}\centerline{\begin{tabular}{|c|c|}\hline
%minimum field containing $H$&number of subgroups\\\hline$H\subseteq k$&$(q-1)(q-2)(q-4)/168$\\
%\hline$H\subseteq k_2$, $H\not\subseteq k$&$q(q-1)(q-2)/8$\\
%\hline$H\subseteq k_3$, $H\not\subseteq k$&$(q-1)(q^2-1)/3$\\
%\hline$H\subseteq k_4$, $H\not\subseteq k_2$&$q^2(q-1)/4$\\\hline\end{tabular}}\be
%Substracting the number of all these polynomials from the total number $q^3-q^2$, we get $|S|=2(q^3-1)/7$.
%This proves the lemma, since all polynomials in $S$ divide the product $(x^{q^3}+x^q+x)(x^{q^3}+x^{q^2}+x)$.
The change of variables, $\tau^n:=x^{q^n}$, $n\ge0$, establishes
an isomorphism between the ring of $k$-linear polynomials in $x$,
with coefficients in $k$ (with composition as the product
operation) and the ring $k[\tau]$ of polynomials in $\tau$. Thus,
the identity:
$$
\tau^7+1=(\tau^3+\tau+1)(\tau^3+\tau^2+1)(\tau+1),
$$ implies that there exists a $k$-linear polynomial $P(x)$ such that $x^{q^7}+x=P(x^{q^3}+x^q+x)$, and
similarly for $x^{q^3}+x^{q^2}+x$. In particular, these two
polynomials divide $x^{q^7}+x$, so that all their irreducible
factors, apart from $x$, have degree seven.
\end{proof}

Let us deal now simultaneously with the cases $\cc_{7,0}$ and
$\cc_{7,1}$. Let $\z\in k_7$ be a root of any fixed polynomial in
$S_0$ (resp. $S_1$). For instance, if $7\nmid [k\colon \ft]$ we
can take $\z$ to be a root of $x^7+x+1$ (resp. $x^7+x^3+1$). We
take,
$$
\bb=\{\ell,\,\ell',\,\ell'',\,\ell''',\,\ell^{iv},\,\ell^v,\,\ell^{vi}\},
$$where $\ell(x,y,z)=\z x+\z'y+\z''z$. We have,
$$
\begin{array}{l}
\ell'''=\ell+\ell',\ \ell^{iv}=\ell'+\ell'',\ \ell^v=\ell+\ell'+\ell'',\ \ell^{vi}=\ell+\ell'',\\
\ell'''=\ell+\ell'',\ \ell^{iv}=\ell+\ell'+\ell'',\
\ell^v=\ell+\ell',\ \ell^{vi}=\ell'+\ell'',
\end{array}
$$respectively in the cases $(7,0)$, $(7,1)$. The rational models will be of the type:
\renewcommand\arraystretch{1.}
$$
N_Q\colon\quad Q^2=\ell\ell'\ell''\ell^v, \qquad \mbox{ resp. }
Q^2=\ell\ell'\ell''\ell^{iv}.
$$
These curves are defined over $k$ if and only if
$Q+Q'=\ell'\ell''$. The seven points of the Fano plane $\bb$ are
conjugate, hence, if $P$ is the intersection point of $\ell$ and
$\ell'$, the set of quadratic forms $Q$ such that $N_Q$ is defined
over $k$ and non-singular is:
$$
\qq_7:=\{Q=(a,b,c,d,e,f)\in (k_7)^6\tq Q+Q'=\ell'\ell'',\,Q(P)\ne
0\}.
$$

For $\ga:=\ga_{7,0}$ (resp.  $\ga:=\ga_{7,1}$), we have
$\aut_k(\bb)=\g_{\t\ga}=\gen{\t\ga}$. One checks easily that
$H_{\t\ga}=\ell'\ell''$, so that
$$
\t\ga^{-1}(Q)=Q^{\t\ga}+\ell'\ell''=(Q^{\t\ga})'.
$$Thus, $\t\ga(Q)=(Q^{\t\ga^{-1}})'$ and the orbit of  $Q=(a,b,c,d,e,f)$ is the cyclic orbit
genera\-ted respectively by:
$$
\t\ga(Q)=(a'+c'+f',a',b',f',d',d'+e'),\quad\
\t\ga(Q)=(b'+c'+e',a',b',d'+f',d',e').
$$
%$$\begin{array}{cccccc}
%a&b&c&d&e&f\\a'+c'+f'&a'&b'&f'&d'&d'+e'\\\Sigma''&a''+c''+f''&a''&d''+e''&f''&d''+f''\\
%b'''+c'''+e'''&\Sigma'''&a'''+c'''+f'''&d'''+f'''&d'''+e'''&d'''+e'''+f'''\\
%a^{iv}+b^{iv}+d^{iv}&b^{iv}+c^{iv}+e^{iv}&\Sigma^{iv}&d^{iv}+e^{iv}+f^{iv}&d^{v}+
%f^{iv}&e^{iv}+f^{iv}\\c^v&a^v+b^v+d^v&b^v+v^v+e^v&e^v+f^v&d^v+e^v+f^v&e^v\\
%b^{vi}&c^{vi}&a^{vi}+b^{vi}+d^{vi}&e^{vi}&e^{vi}+f^{vi}&d^{vi},\end{array}$$ in the case
%(7,0) and in the case (7,1):
%$$\begin{array}{cccccc}
%a&b&c&d&e&f\\b'+c'+e'&a'&b'&d'+f'&d'&e'\\a''+b''+d''&b''+c''+e''&a''&d''+e''+f''&d''+f''&d''\\
%\Sigma'''&a'''+b'''+d'''&b'''+c'''+e'''&e'''+f'''&d'''+e'''+f'''&d'''+f'''\\
%a^{iv}+c^{iv}+f^{iv}&\Sigma^{iv}&a^{iv}+b^{iv}+d^{iv}&d^{iv}+e^{iv}&e^{v}+
%f^{iv}&d^{iv}+e^{iv}+f^{iv}\\c^v&a^v+c^v+f^v&\Sigma^v&f^v&d^v+e^v&e^v+f^v\\
%b^{vi}&c^{vi}&a^{vi}+c^{vi}+f^{vi}&e^{vi}&f^{vi}&d^{vi}+e^{vi}.
%\end{array}$$

We have $\aut_k(N_Q)=\{1\}$ for all $N_Q$ except for one, which is
a twist of the Klein quartic (this was seen in paragraph \ref{nbr}
too). This special model $N_Q$ corresponds to:
$$
Q=\ell^2+(\ell')^2+(\ell'')^2+\ell\ell'+\ell'\ell''+\ell\ell'',
$$and it has $\aut(N_Q)=\gen{\t\ga}\simeq C_7$.

\section{Non-ordinary curves}\label{sec:2}
As in the ordinary case, our starting point is a family of models
representing all $\kb$-isomorphy classes of non-ordinary curves
\cite{wa} \cite{r}.

\begin{prop}\label{dosu}
Let $C\subseteq \pr{2}$ be a non-ordinary, non-singular quartic
plane curve defined over $\kb$. Then $C$ satisfies one of the
following equivalent conditions:
\begin{enumerate}
\item The Jacobian $J_C$ of $C$ has 2-rank 2, resp. 1, resp. 0.
\item $C$ has 4, resp. 2, resp. 1 bitangents.
\item $C$ is isomorphic to a curve $C_Q$ with equation:
$$
Q^2=xyz(y+z);\quad \mbox{ resp. }\ Q^2=xy(y^2+xz);\quad  \mbox{
resp. }\ Q^2=x(y^3+x^2z),
$$where $Q(x,y,z)=ax^2+by^2+cz^2+dxy+eyz+fzx$ is a quadratic form satisfying:
$$
abc\ne0,\, b+c+e\ne0;\quad \mbox{ resp. }\  ac\ne0;\quad \mbox{
resp. }\  c\ne 0.\qed
$$
\end{enumerate}
\end{prop}

In this section we find several families $\nn$ of rational normal
models containing disjoint sets of $k$-isomorphism classes, such
that every non-ordinary quartic is $k$-isomorphic to some normal
model. For each family $\nn$ we find a finite group $G$ acting
fully by $k$-isomorphisms on $\nn$; this allows us to explicitely
determine the $k$-isomorphy classes and to compute the
$k$-automorphism group of each curve. Also, we obtain the number
$|G\backslash \nn|$ of $k$-isomorphism classes contained in $\nn$
by applying (\ref{comb}).

As before, we abuse of language and denote the quadratic forms
(\ref{qf}) in $x,y,z$ simply by $Q=(a,b,c,d,e,f)$.

\subsection{Curves with four bitangents}
Let $\qq$ be the set of all quadratic forms $Q=(a,b,c,d,e,f)$
defined over $\kb$
 such that $abc\ne0$, $b+c+e\ne0$. Let $\cc_{\qq}$ be the family of all curves:
$$
C_Q\colon\quad Q(x,y,z)^2=xyz(y+z),\qquad Q\in\qq.
$$
The set of bitangents of these curves is $\bb=\{x,y,z,y+z\}$. As
in the ordinary case, the group $\aut(\bb)$ acts fully by
isomorphisms on the family $\cq$; in fact, the elements of
$\aut(\bb)$ are the only transformations that preserve the quartic
$xyz(y+z)$ modulo squares. This group is isomorphic to $S_3\times
{\kb}^*$ \cite{wa}, where the subgroup isomorphic to $S_3$
contains the automorphisms permuting $y,\,z,\,y+z$:
\begin{multline*}
\{1,\,\ga_{\tau}(x,y,z)=(x,y,y+z),\,\ga_{\tau_2}(x,y,z)=(x,z,y),\,\ga_{\tau_3}(x,y,z)=(x,y+z,z),\\
 \ga_{\rho}(x,y,z)=(x,z,y+z),\,\ga_{\rho^2}(x,y,z)=(x,y+z,y)\},
\end{multline*}
and the subgroup isomorphic to $\kb^*$ consists of the
transformations:
$$
\ga_t(x,y,z)=(t^3x,t^{-1}y,t^{-1}z),\qquad t\in \kb^*.
$$
Every element of $\aut(\bb)$ can be written in a unique way as
$\ga_{\beta,t}:=\ga_t\circ\ga_{\beta} =\ga_{\beta}\circ \ga_t$,
with $\beta\in S_3$, $t\in \kb^*$. Moreover,
$$\ga_{\beta',t'}\circ \ga_{\beta,t}=\ga_{\beta'\beta,t't},\quad\forall \beta,\beta'\in S_3,\,
\forall t,t'\in \kb^*.$$

For any $\beta\in S_3$ let us denote by $\dd_{\beta}$ the set of
all descent data over $k$ of curves of $\cc_{\qq}$ generated by
$\ga_{\beta}$ (cf. paragraph \ref{ddata}):
$$\dd_{\beta}=\{(C,\ga_{\beta})\tq C\in\cc_{\qq},\,\ga_{\beta}(C)=\si C\}.$$

\begin{prop}\label{dquatre}
$$
\dsc(\cc_{\qq})=\dsc(\dd_1)\amalg\dsc(\dd_{\tau})\amalg\dsc(\dd_{\rho}).
$$
Moreover, two descent data in the same family $\dd_{\beta}$ are
equivalent iff they are in the same orbit under the action of the
centralizer of $\ga_{\beta}$ in $\aut_k(\bb)\simeq S_3\times k^*$.
\end{prop}

\begin{proof}
Assume that $\ga_{\beta,t}\colon C_Q\iso C_{\si Q}$ generates
descent data over $k$ of degree $n$ for the curve $C_Q$. Clearly:
$$
\sij{n-1}\ga_{\beta,t}\circ \cdots\circ\si \ga_{\beta,t}\circ
\ga_{\beta,t}=1 \imp \sij{n-1}t\cdots\si tt=1.
$$ By Hilbert's 90th theorem, there exists $s\in \kb^*$ such that $t=s/\si s$. In particular,
$\si\ga_{1,s}\circ\ga_{\beta,t}\circ\ga_{1,s}^{-1}=\ga_{\beta,1}$.
Thus, by (\ref{eqdd}), our original descent datum
$(C_Q,\ga_{\beta,t})$ is equivalent to some descent datum having
$t=1$.

Once restricted to descent data with $t=1$, the situation is
completely analogous to that of paragraph \ref{ddata}. The
assertion of the proposition derives from results analogous to
Proposition \ref{prop3}
 and Theorem \ref{bij}.
\end{proof}

Let $k_2$, $k_3$ denote respectively the quadratic and cubic
extension of $k$ in $\kb$. We fix $u\in k_2\setminus k$ with
equation $u^2+u=r$, for some $r\not\in\as(k)$. Also, we fix $v\in
k_3\setminus k$ with equation $v^3+v=s$ for adequate $s\in k^*$.
By\cite[Lemma 7]{cnp} we have $v'=s^{-1}v^2+tv$,
$v''=s^{-1}v^2+(t+1)v$, for certain $t\in k$ satisfying
$t^2+t+1=s^{-1}$.

By Proposition \ref{dquatre} there is still an action of $k^*$ on
all descents of the family $\cq$, given by the isomorphisms
$\ga_{1,t}$, $t\in k^*$. Thus, we can reduce the family of normal
models by assuming that the coefficient $a$ in the quadratic forms
runs on a fixed system of representatives of $k^*/(k^*)^3$. So, we
take the following sets of quadratic forms:
\renewcommand\arraystretch{1.4}
$$
\begin{array}{l}
\tilde\qq_1:=\{(a,b,c,d,e,f)\in k^6\tq abc\ne0,\, b+c+e\ne 0\},\\
\tilde\qq_2:=\{(a,b,c,d,e,f)\in k^6\tq a\ne0,\, (b,e)\ne(cr,c)\},\\
\tilde\qq_3:=\{(a,b,c,d,e,f)\in k^6\tq a\ne0,\, (b,c,e)\ne
(0,0,0)\},
\end{array}
$$
and subsets $\qq_i:=\{Q\in\tilde\qq_i\tq a\in k^*/(k^*)^3\}$.
Consider the families of rational normal models:
$$\tilde\nn^{(4)}_i:= \{N_Q\tq Q\in\tilde\qq_i\},\quad
\nn^{(4)}_i:=\{N_Q\tq Q\in\qq_i\}, \qquad i=1,2,3,$$ where $N_Q$
is respectively:
$$
\begin{array}{ll}
N_Q\colon \quad Q^2=xyz(y+z),&\\
N_Q\colon \quad Q^2=xy\ell\ell'=xy(ry^2+yz+z^2), &\ell(x,y,z)=uy+z,\\
N_Q\colon \quad Q^2=x\ell\ell'\ell''=x(y^3+ty^2z+(t+1)yz^2+z^3),
 &\ell(x,y,z)=y+v'v^{-1}z,
\end{array}
$$

\renewcommand\arraystretch{1.}
The isomorphism:
$$
\phi(x,y,z)=(x,y,z),\ \mbox{ resp. }(x,y,\ell(x,y,z)) \mbox{ resp.
}(s^{-1}x,v\ell(x,y,z),v'\ell'(x,y,z))
$$sets the families of models $\tilde\nn^{(4)}_i$, $i=$1,2,3 in respective bijection
with the descent data $\dd_{\beta}$, for $\beta=1,\tau,\rho$. By
Proposition \ref{dquatre} and a result analogous to Lemma
\ref{aut}, the group $\,\phi^{-1}\circ \left((S_3)_{\beta}\times
k^*\right)\circ \phi\,$ acts fully by $k$-isomorphisms on the
corresponding family $\tilde\nn^{(4)}_i$. Hence, the group
$\,G=\phi^{-1}\circ \left((S_3)_ {\beta}\times
\mu_3(k)\right)\circ \phi\,$ acts fully by $k$-isomorphisms on the
corresponding family $\nn^{(4)}_i$. Here and in the sequel,
$\mu_n(k)$ denotes the group of $n$-th roots of unity that are
contained in $k$.

By Proposition \ref{dquatre}, every non-singular quartic $C$ over
$k$ with four bitangents is $k$-isomorphic to a curve in one of
the three families $\nn^{(4)}_i$. The minimum field of definition
of the set of bitangents tells us to which family belongs the
$k$-isomorphism class of $C$. Any transformation in $\pg{k}$
sending the set of four bitangents of $C$ to the set of four lines
that has been chosen for the family $\nn^{(4)}_i$, will take $C$
into one of the normal models $N_Q$.

We describe now in each case the $k$-isomorphism classes and the
$k$-automorphism group. In the sequel we denote
$N:=|\mu_3(k)|=|k^*/(k^*)^3|$; that is, $N=3,1$, according to $q$
being a square
 or not.

\subsubsection*{Normal models $\nn^{(4)}_1$}
Since $\phi=1$, we have $G=\{\ga_{\beta,t}\tq\beta\in
S_3,\,t^3=1\}\simeq S_3\times \mu_3(k)$. For any given
$Q=(a,b,c,d,e,f)\in\qq_1$ and $t\in\mu_3(k)$, we have
$$
\ga_{1,t}(a,b,c,d,e,f)=(a,t^2b,t^2c,td,t^2e,tf),
$$which has no fixed points if $t\ne 1$. On the other hand,
the orbit of $Q$ under the action of $S_3$ is:
$$
\begin{array}{cccc}
(a,b,c,d,e,f)&\quad (1)&\qquad(a,B,c,d+f,e,f)&\quad (\ga_{\tau})\\
(a,c,b,f,e,d)&\quad (\ga_{\tau_2})&\qquad(a,b,B,d,e,d+f)&\quad (\ga_{\tau_3})\\
(a,B,b,d+f,e,d)&\quad (\ga_{\rho})&\qquad(a,c,B,f,e,d+f)&\quad
(\ga_{\rho^2}),
\end{array}
$$where we have denoted $B:=b+c+e$. This leads immediately to:
$$
\aut_k(N_Q)=\left\{
\begin{array}{ll}
\{\ga_{\beta,1}\}\simeq S_3,&\mbox{ if }b=c=e,\,d=f=0,\\
\gen{\tau}\simeq C_2,&\mbox{ if }c=e,\,f=0\ne d,\\
\gen{\tau_2}\simeq C_2,&\mbox{ if }b=c,\,d=f\ne0,\\
\gen{\tau_3}\simeq C_2,&\mbox{ if }b=e,\,d=0\ne f,\\
\{1\},&\mbox{ otherwise}.
\end{array}
\right.
$$

Finally, we have $|\qq_1|=N(q-1)^3q^2$ and\renewcommand\arraystretch{1.4}
$$
\begin{array}{ll}
\qq_1(\ga_{\tau})=\{(a,b,c,d,c,0)\in k^6\tq a\in k^*/(k^*)^3,\,bc\ne0\}&\imp|\qq_1(\tau)|=N(q-1)^2q,\\
\qq_1(\ga_{\rho})=\{(a,c,c,0,c,0)\in k^6\tq a\in
k^*/(k^*)^3,\,c\ne0\}, &\imp|\qq_1(\ga_{\rho})|=N(q-1).
\end{array}\renewcommand\arraystretch{1.}
$$
We get from (\ref{comb}):
\begin{multline*}
|G\backslash \nn^{(4)}_1|=\frac 1{6N}\left(N(q-1)^3q^2+3N(q-1)^2q+2N(q-1)\right)=\\
=\frac 16\left(q^5-3q^4+6q^3-7q^2+5q-2\right).
\end{multline*}

\subsubsection*{Normal models $\nn^{(4)}_2$}
Since $\phi^{-1}\circ \ga_{\tau}\circ\phi=\ga_{\tau}$, we have
$G=\{\ga_{\beta\,t}\tq \beta\in\gen{\tau}, \,t^3=1\}\simeq
C_2\times \mu_3(k)$. Since the elements of $G$ leave the quartic
$xy\ell\ell'$ invariant, the action on $\qq_2$ is the same than
before. Hence, all $k$-automorphism groups are trivial except for:
$$
\aut_k(N_Q)=\gen{\ga_{\tau}}\simeq C_2,\quad \mbox{ if }c=e,\,f=0.
$$

Finally, we have $|\qq_2|=N(q-1)(q^2-1)q^2$,
$|\qq_2(\ga_{\tau})|=N(q-1)^2q$ ($a\ne0$, $b\ne cr$, $d$
arbitrary) and we get from (\ref{comb}):
$$
|G\backslash \nn^{(4)}_2|=\frac
1{2N}\left(Nq(q-1)^2(q^2+q+1)\right)=\frac
12\left(q^5-q^4-q^2+q\right).
$$

\subsubsection*{Normal models $\nn^{(4)}_3$}
Since $\phi$ is symmetric, we have $\phi^{-1}\circ
\ga_{\rho}\circ\phi=\t\ga_{\rho}=\ga_{\rho}$. Hence,
$G=\gen{\ga_{\rho}}\times \mu_3(k)\simeq C_3\times \mu_3(k)$.
Since the elements of $G$ leave the quartic $x\ell\ell'\ell''$
invariant, the action on $\qq_3$ is the one described in the case
$\nn^{(4)}_1$. Hence, all $k$-automorphism groups are trivial
except for:
$$
\aut_k(N_Q)=\gen{\ga_{\rho}}\simeq C_3,\quad\mbox{ if
}b=c=e,\,d=f=0.
$$
Finally, we have $|\qq_3|=N(q^3-1)q^2$,
$|\qq_3(\ga_{\rho})|=N(q-1)$ and we get from (\ref{comb}):
$$
|G\backslash \nn^{(4)}_3|=\frac
1{3N}\left(N(q^3-1)q^2+2N(q-1))\right)=\frac
13\left(q^5-q^2+2q-2\right).
$$

\subsection{Curves with two bitangents}
Let $\cq$ be the family of non-singular quartics:
$$
C_Q\colon\quad Q^2=xy(y^2+xz),
$$ where $Q(x,y,z)=ax^2+by^2+cz^2+dxy+eyz+fzx$ is a quadratic form with
coefficients in $\kb$ satisfying $ac\ne0$.

The subgroup of automorphisms of $\pr{2}$ that preserve the
quartic $xy(y^2+xz)$ modulo squares consists of the
transformations \cite{wa}:
$$
\ga_{t,u}(x,y,z)=(t^3x,t^{-1}y,t^{-5}(z+u y)),\qquad
t\in\kb^*,\,u\in\kb.
$$The composition of two such transformations is given by:
$$
\ga_{t',u'}\circ \ga_{t,u}=\ga_{t't,u+t^4u'},\qquad \forall
t,t'\in\kb^*,\,u,u'\in \kb.
$$Clearly, this group acts fully by isomorphisms on the family $\cq$.

\begin{prop}\label{trdes}
Every non-singular quartic plane curve $C\subseteq\pr{2}$ with two
bitangents is $k$-isomorphic to a curve $C_Q$ for some $Q$ with
coefficients in $k$.
\end{prop}

\begin{proof}
The assertion is true over $\kb$ by Proposition \ref{dosu}. We
need only to check that all descent data over $k$ of these curves
$C_Q$ are trivial. Assume that $\ga_{t,u}\colon C_Q\iso C_{\si Q}$
generates descent data of degree $n$ for the curve $C_Q$ defined
over $\kb$. Arguing as in the proof of Proposition \ref{dquatre}
we can assume that $t=1$. Then,
 taking $v\in\kb$ such that $v+\si v=u$ we have: $\si\ga_{1,v}\circ\ga_{1,u}
\circ\ga_{1,v}^{-1}=1$, so that all these descent data are
trivial.
\end{proof}

On the family of curves $Q^2=xy(y^2+xz)$ with $Q$ defined over $k$
we can still apply $k$-isomorphisms $\ga_{t,u}$ with $t\in k^*$,
$u\in k$. We have, $\ga_{t,u}(C_Q)=C_{\ga_{t,u}(Q)}$, where, for
$Q=(a,b,c,d,e,f)$:
\begin{equation}\label{orbit}
\ga_{t,u}(Q)=(at^{-6},(b+eu+cu^2)t^2,ct^{10},(d+\sqrt
u+fu)t^{-2},et^6,ft^2).
\end{equation}
This allows us to consider a smaller family of curves still
containing all $k$-isomorphy classes. By choosing a suitable $t$
we can prefix the class of $a$ modulo cubes. By choosing a
suitable $u$ we can assume that $d=0$ (if $f=0$ or $df\in\as(k)$)
or $df=r_0$, for some fixed element $r_0\in k\setminus\as(k)$.
This leads to the consideration of two families of rational normal
models:
\renewcommand\arraystretch{1.4}
$$
\begin{array}{l}
\qq_1:=\{Q=(a,b,c,d_0,e,f)\in k^6\tq a\in k^*/(k^*)^3,\,c\ne0,\,d_0\in\{0,f^{-1}r_0\}\},\\
\qq_0:=\{Q=(a,b,c,0,e,0)\in k^6\tq a\in k^*/(k^*)^3,\,c\ne0\},\\
\nn_1^{(2)}:=\{C_Q\tq Q\in\qq_1\},\qquad \nn^{(2)}_0:=\{C_Q\tq
Q\in\qq_0\}.
\end{array}
$$
\renewcommand\arraystretch{1.4}

\subsubsection*{Normal models $\nn_1^{(2)}$}
The group $\,G:=\{\ga_{t,u}\tq
t^3=1,\,u\in\{0,f^{-2}\}\}\simeq\mu_3(k)\times C_2\,$ acts fully
by $k$-isomorphisms on $\nn_1^{(2)}$. From (\ref{orbit}) we have:
$$
\ga_{t,u}(a,b,c,d_0,e,f)=(a,(b+eu+cu^2)t^2,ct,d_0t^{-2},e,ft^2).
$$Therefore, all $k$-automorphisms of the curves in $\nn_1^{(2)}$ are trivial except for:
$$\aut_k(N_Q)=\gen{\ga_{1,f^{-2}}}\simeq C_2,\quad\mbox{ if }e=cf^{-2}.
$$

Clearly, $|\nn_1^{(2)}|=2N(q-1)^2q^2$, the element
$\ga_{1,f^{-2}}$ has $2N(q-1)^2q$ fixed points and all other
transformations in $G$ have no fixed points. By (\ref{comb}):
$$
|G\backslash \nn_1^{(2)}|=\frac
1{2N}\left(2N(q-1)^2q^2+2N(q-1)^2q\right)=(q-1)^2q(q+1).
$$

\subsubsection*{Normal models $\nn^{(2)}_0$}
The group $G:=\{\ga_{t,0}\tq t^3=1\}\simeq\mu_3(k)$ acts fully by
$k$-isomorphisms on $\nn^{(2)}_0$. From (\ref{orbit}) we have now:
$$
\ga_{t,0}(a,b,c,0,e,0)=(a,bt^2,ct,0,e,0).
$$Therefore, all $k$-automorphism groups are trivial and
$$
|G\backslash \nn^{(2)}_0|=\frac 1{|G|}|\nn^{(2)}_0|=(q-1)q^2.
$$

\subsection{Curves with one bitangent}
Let $\cq$ be the family of non-singular quartics:
$$
C_Q\colon\quad Q^2=x(y^3+x^2z),
$$ where $Q(x,y,z)=ax^2+by^2+cz^2+dxy+eyz+fzx$ is a quadratic form with
coefficients in $\kb$ satisfying $c\ne0$.

On $\cq$ we have a full action by isomorphisms by the subgroup of
automorphisms of $\pr{2}$ that preserve the quartic $x(y^3+x^2z)$
modulo squares. This group consists of the transformations
\cite{wa}:
$$
\ga_{t,u,v}(x,y,z)=(t^3x,t^{-1}(y+ux),t^{-9}(z+u^2y+vx)),\qquad
t\in\kb^*,\,u,v\in\kb.
$$The composition of two such transformations is given by:
$$
\ga_{t',u',v'}\circ
\ga_{t,u,v}=\ga_{t't,u+t^4u',v+u(u')^2t^8+v't^{12}},\qquad \forall
t,t'\in\kb^*, \,u,v,u',v'\in \kb.
$$

\begin{prop}\label{trdes2}
Every non-singular quartic plane curve $C\subseteq\pr{2}$ with one
bitangent is $k$-isomorphic to a curve $C_Q$ for some $Q$ with
coefficients in $k$.
\end{prop}

\begin{proof}
As in Proposition \ref{trdes} we need only to check that all
descent data over $k$ of these curves $C_Q$ are trivial.

Assume that $\ga_{t,u,v}\colon C_Q\iso C_{\si Q}$ generates
descent data for the curve $C_Q$ defined over $\kb$. By the
argument we used in the proof of Proposition \ref{trdes} we can
assume that $t=1$, $u=0$. Again, these descent data are all
trivial, since taking $w\in\kb$ such that $w+\si w=v$ we have:
$\si\ga_{1,0,w}\circ\ga_{1,0,v}\circ\ga_{1,0,w}^{-1}=1$.
\end{proof}

On the family of curves $Q^2=x(y^3+x^2z)$ with $Q$ defined over
$k$ we can still apply $k$-isomorphisms $\ga_{t,u,v}$ with $t\in
k^*$, $u,v\in k$. We have, $\ga_{t,u,v}(C_Q)=C_{\ga_{t,u,v}(Q)}$,
where, for $Q=(a,b,c,d,e,f)$:
\begin{multline}\label{orbit2}
\ga_{t,u,v}(Q)=((a+bu^2+c(v+u^3)^2+du+(eu+f)(v+u^3)+\sqrt v)t^{-6},\\
(b+eu^2+cu^4)t^2,ct^{18},(d+ev+fu^2+\sqrt
u)t^{-2},et^{10},(eu+f)t^6).
\end{multline}
This allows us to consider a smaller family of curves still
containing all $k$-isomorphy classes. By choosing a suitable $t$
we can prefix the class of $c$ modulo $(k^*)^9$. If $e\ne0$, we
can achieve $d=f=0$ by choosing suitable $u,v$; if $e=0$ we obtain
$b=0$ by a suitable choice of $u$.  This leads to the
consideration of two families of rational normal models. The
family $\nn_1^{(1)}$ consists of all curves with equation:
$$
N_Q\colon\qquad Q^2=x(y^3+x^2z),\qquad Q=(a,b,c,0,e,0),\ c\in
k^*/(k^*)^9,\,a,b\in k,\,e\in k^*,
$$and the family $\ss:=\nn_0^{(1)}$ gathers all curves with equation:
$$
N_Q\colon\qquad Q^2=x(y^3+x^2z),\qquad Q=(a,0,c,d,0,f),\ c\in
k^*/(k^*)^9,\,a,d,f\in k.
$$
We shall see below that a quartic defined over $k$ is
supersingular if and only if it is $k$-isomorphic to a curve in
the family $\ss$.

>From now on, we let $N:=|\mu_9(k)|=|k^*/(k^*)^9|$; that is,
$N=9,3,1$, according respectively to $q\equiv 1\md9$, $q\equiv
4,7\md9$ or $q\not\equiv 1\md 3$.

\subsubsection*{Normal models $\nn_1^{(1)}$}
The group $G:=\{\ga_{t,0,0}\tq t^9=1\}\simeq\mu_9(k)$  acts fully
by $k$-isomorphisms on $\nn_1^{(1)}$. From (\ref{orbit2}) we have:
$$
\ga_{t,0,0}(a,b,c,0,e,0)=(at^3,bt^2,c,0,et,0).
$$Therefore, the $k$-automorphism groups are all trivial and
$$
|G\backslash \nn_1^{(1)}|=\frac 1N|\nn_1^{(1)}|=(q-1)q^2.
$$

\subsubsection*{Normal models for supersingular quartics}
The group $G:=\{\ga_{t,0,v}\tq t^9=1,\,v\in k\}\simeq
\mu_9(k)\times k$  acts fully by $k$-isomorphisms on $\ss$. From
(\ref{orbit2}) we have now:
$$
\ga_{t,0,v}(a,0,c,d,0,f)=((a+cv^2+fv+\sqrt
v)t^3,0,c_0,dt^{-2},0,ft^{-3}).
$$
We could reduce the family of normal models by choosing adequate
representatives for the coefficient $a$ modulo this action, but
the analysis of the case is easier if we don't do this. For fixed
$c\in k^*$, $f\in k$, let us consider the $\ft$-linear
homomorphism,
$$E_{c,f}\colon k\to k, \qquad E_{c,f}(x)=cx^2+fx+\sqrt x.$$ One checks immediately
that:\renewcommand\arraystretch{1.4}
$$
\aut_k(N_Q)\!=\!\left\{\negmedspace\negmedspace
\begin{array}{l}
\ker(E_{c,f}), \ \qquad\qquad\qquad\qquad\qquad\qquad\qquad\qquad\qquad\mbox{ if }d\ne0,\\
\mathbb{\mu}_3(k)\times \ker(E_{c,f}),\ \
\,\qquad\qquad\qquad\qquad\qquad
\qquad\qquad\mbox{ if }d=0,\,f\ne0,\\
\left(\mu_3(k)\times\ker(E_{c,0})\right)\amalg
\{\ga_{t,0,v}\!\in\! G\tq t^3\ne1,\,v\in
E_{c,0}^{-1}(t^3a)\},\,\mbox{if }d=f=0,
\end{array}\right.
$$
\renewcommand\arraystretch{1.}
For fixed $c\in k^*$, $v\in k$, denote  $$\nv:=|\{f\in k\tq v\in
\ker(E_{c,f})\}|= \left\{\begin{array}{ll}q,&\qquad \mbox{ if
}v=0,\\1,&\qquad\mbox{ if }v\ne0,\end{array}\right.
$$which is independent of $c$. Clearly $\sum_{v\in k}\nv=2q-1$.
The number of fixed points of any $\ga_{t,0,v}\in\g$ is:
$$
\begin{array}{cll}
N\,q^2\,\nv,&\qquad\mbox{ if }t=1,&\qquad (a,d\in k,\,c\in k^*/(k^*)^9))\\
N\,q\,\nv,&\qquad\mbox{ if }t^3=1,\,t\ne1,&\qquad (d=0,\,a\in k,\,c\in k^*/(k^*)^9))\\
N,&\qquad\mbox{ if }t^3\ne1.&\qquad
(f=d=0,\,a=t^{-3}E_{c,0}(v),\,c\in k^*/(k^*)^9)
\end{array}
$$
We can apply now (\ref{comb}) to count the number of $k$-isomorphy
classes:
$$
|G\backslash \ss|=\frac 1{Nq}\sum_{t,v}|\ss(\ga_{t,0,v})|= \frac
1{Nq}\sum_{v\in k}N\left((q^2+\delta q)\nv+
\eps\right)=(q+\delta)(2q-1)+\eps,
$$where $\delta=|\mu_3(k)|-1$ and $\eps=|\mu_9(k)|-|\mu_3(k)|$.

\subsubsection*{Stratification by the Newton polygon} For an
abelian variety $A$ over a finite field $\fq$
%, the Newton polygon (NP) coincides with the NP of the characteristic
%polynomial of the Frobenius endomorphism of $A$, with respect to
%the $2$-adic valuation $v_2$ normalized by $v_2(q)=1$.
let us denote by $\op{Slp}(A)$ the set of slopes of the Newton
polygon (NP) of $A$ \cite{oo}. It is well known that
\renewcommand\arraystretch{1.2}
\begin{equation}\label{slopes}
\begin{array}{l}
\op{Slp}(A_1\times A_2)=\op{Slp}(A_1)\cup \op{Slp}(A_2), \\
\dim(A)\le 2\imp \op{Slp}(A)\subseteq\{0,1/2,1\}.
\end{array}
\end{equation}
\renewcommand\arraystretch{1.}

If $\dim(A)=3$, the NP is determined by the 2-rank of $A$, except
for the abelian varieties of 2-rank zero, whose NP has either two
sides with slopes 1/3, 2/3 or one side with slope 1/2. We say
respectively that $A$ is ``of type 1/3" or ``supersingular". We
use the same terminology for a quartic with one bitangent
according to the shape of the NP of its jacobian.

\begin{prop}\label{ss}
\casos\cas All quartics in the family $\ss$ are supersingular.\cas
All quartics in the family $\nn_1^{(1)}$ are of type 1/3. \fcasos
\end{prop}

\begin{proof}
a)  Let $C=N_Q$ be a quartic in the family $\ss$. By taking $x=0$
as the line at infinity we obtain an affine model:
$$
C\colon\quad cz^4+fz^2+z=y^3+dy^2+a.
$$
Let $v\in\kb^*$ be a non-trivial root of the polynomial
$cz^4+fz^2+z$. The curve $C$ admits the involution
$\varphi(x,y,z)=(x,y,z+v)$ and the quotient curve $C':=C/\varphi$
has genus one. In fact, working with the affine model, the
morphism $C\to C'$ can be seen at the level of function fields as:
$$
k(y,u)\subseteq k(y,z),\qquad u:=z(z+v).
$$Since
$cu^2+v^{-1}u=cz^4+fz^2+z$, we can take
$$
cu^2+v^{-1}u=y^3+dy^2+a,
$$ as an affine model for $C'$ and this is clearly an elliptic
curve.

Therefore, the Jacobian of $C$ admits a non-trivial morphism to an
elliptic curve and it cannot be of type 1/3 because an abelian
threefold of type 1/3 is absolutely simple by (\ref{slopes}).\e

b) Let $\mm_3\hookrightarrow \aa_{3,1}$ be the respective moduli
spaces of curves of genus three and principally polarized abelian
threefolds over $\ft$. The supersingular locus in $\aa_{3,1}$ is a
closed subset \cite[(2.4)]{oo} and T. Katsura and F. Oort have
proved that it has dimension 2 and is absolutely irreducible
\cite[section 6]{ko}.

Let $SS$ be the supersingular locus of $\mm_3$ and let $S\subseteq
\mm_3$ be the set of images in $\mm_3$ of all curves of the family
$\ss$ for all finite fields of characteristic 2. By a) we have
$S\subseteq SS$. Since all curves in $\ss$ have non-trivial
automorphisms, they lie in the singular locus of the moduli space:
$S\subseteq (\mm_3)_{\mbox{\scriptsize $\op{sing}$}}$ .

By Lemma \ref{moduli} below, we have $|S(\fq)|=q^2$ for all $q$.
Hence, the result of Katsura-Oort implies that $\overline{S}=SS$.
Since the singular locus of a variety is a closed subset, we get
$SS=\overline{S}\subseteq (\mm_3)_{\mbox{\scriptsize
$\op{sing}$}}$. On the other hand, all curves in $\nn_1^{(1)}$
furnish smooth points of $\mm_3$ since their automorphism groups
are trivial; hence, none of these curves is supersingular.
\end{proof}

\section{Number of curves of genus three}\label{sec:3}
Gathering all computations of section 1 and section 2 we have:

\begin{teor}\label{total} There are
$$q^6+q^4-q^3+2q^2+4-[4q-2]_{_{q\equiv-1\,(\op{mod}\,3)}}+
[6]_{_{q\equiv1\,(\op{mod}\,9)}}$$ $\fq$-isomorphism classes of
non-singular quartic plane curves defined over $\fq$. According to
the Newton polygon of the Jacobian, these curves are distributed
in the following way:\be

\renewcommand\arraystretch{1.4}
\centerline{\begin{tabular}{|c|c|} \hline
ordinary&$q^6-q^5+q^4-3q^3+5q^2-6q+7$\\
\hline 2-rank two&$q^5-q^4+q^3-2q^2+2q-1$\\\hline
2-rank one&$q^4-2q^2+q$\\
\hline type 1/3&$q^3-q^2$\\
 \hline supersingular&$2q^2-q
+[4q-2]_{_{q\equiv1\,(\op{mod}\,3)}}+[6]_{_{q\equiv1\,(\op{mod}\,9)}}$\\
\hline
\end{tabular}}
\end{teor}\be
\renewcommand\arraystretch{1.}

The notation $[a]_{\mbox{\scriptsize condition}}$ indicates that
$a$ has to be added in the formula if ``{\it condition}'' is
satisfied.

The Newton polygon of a hyperelliptic curve of genus 3 is
completely determined by the number of Weierstrass points. On one
hand, for any hyperelliptic curve $C$ defined over a perfect field
of characteristic 2 we have:
$$
\dim_{\ft}(J_C[2])=|\ww|-1,
$$where $J_C$ is the Jacobian of $C$ and $\ww$ the set of
Weierstrass points. On the other hand, J. Scholten and H.J. Zhu
have proved that there are no supersingular hyperelliptic curves
of genus 3 in characteristic 2 \cite{scholten}. Thus, when $C$ has
2-rank zero it is necessarily of type 1/3.

In \cite[Table 2]{ns} there are formulas for the number of
hyperelliptic curves of genus 3 with prescribed ramification
divisor. Adding these computations to the results of Theorem
\ref{total} we get:

\begin{cor}
There are
$$q^6+2q^5+q^4+q^3+q^2+q+2-[4q-2]_{_{q\equiv-1\,(\op{mod}\,3)}}+
[6]_{_{q\equiv1\,(\op{mod}\,9)}}+[12]
_{_{q\equiv1\,(\op{mod}\,7)}}$$ $\fq$-isomorphism classes of
smooth projective curves of genus three defined over $\fq$.
According to the Newton polygon of the Jacobian, these curves are
distributed in the following way:\be

\renewcommand\arraystretch{1.4}
\centerline{\begin{tabular}{|c|c|} \hline
ordinary&$q^6+q^5-q^4-q^3+q^2-4q+7$\\
\hline 2-rank two&$q^5+q^4-3q^3+q^2+q-1$\\\hline
2-rank one&$q^4+4q^3-4q^2+q-2$\\
\hline type 1/3&$q^3+q^2+[12]
_{_{q\equiv1\,(\op{mod}\,7)}}$\\
 \hline supersingular&$2q^2-q
+[4q-2]_{_{q\equiv1\,(\op{mod}\,3)}}+[6]_{_{q\equiv1\,(\op{mod}\,9)}}$\\
\hline
\end{tabular}}
\end{cor}\be
\renewcommand\arraystretch{1.}

Also, from the computations of sections 1,2 we can deduce the
number of $k$-rational points in the moduli space. Let us quote
first a result that is an immediate consequence of
$\cite[5.1]{vd}$.

\begin{lem}\label{moduli}
Let $\cc$ be a finite family of projective smooth curves of genus
$g>1$, defined over a finite field $k$ with the property that
every $k$-curve of genus $g$ that is isomorphic to some curve in
$\cc$ is $k$-isomorphic to some curve in $\cc$. Suppose that $\cc$
is the disjoint union, $\cc=\amalg_i\cc_i$, of subfamilies $\cc_i$
admitting a full action by $k$-isomorphisms by a finite group
$G_i$. Then, the image of the family $\cc$ in the subset of
$k$-rational points of the moduli space of curves of genus $g$ has
cardinality $\sum_i|\cc_i|/|G_i|$.
\end{lem}

\begin{proof}
By \cite[5.1]{vd}, the contribution of the curves in this family
to the moduli space is the weighted sum
$\sum_{C}|\aut_k(C)|^{-1}$, for $C$ running on a system of
representatives of $k$-isomorphism classes of the curves in $\cc$.
If we have one single group $G$ acting fully by $k$-isomorphisms
on the family $\cc$, this weighted sum can be computed as:
$$\sum_{C\in G\backslash\cc}|\aut_k(C)|^{-1}=\sum_{C\in G\backslash\cc}|G(C)|^{-1}=
\sum_{C\in G\backslash\cc}\frac {|G\cdot C|}{|G|}=\frac
{|\cc|}{|G|},
$$where $G(C)$ is the stabilizer of $C$ and $G\cdot C$ the orbit of $C$. In the general case, the
contribution to the moduli space would be $\sum_i\sum_{C\in
G_i\backslash\cc_i}|\aut_k(C)|^{-1}$ and the above computation
applies to  each subfamily $\cc_i$.
\end{proof}

This result can be applied to each of the five families obtained
by gathering all normal models with the same Newton polygon, to
get closed formulas for the cardinalities of the strata by the
Newton polygon of the non hyperelliptic locus
$\mm^{\mbox{\scriptsize $\op{nh}$}}_3$ of the moduli space. Adding
to these computations the results of \cite[Table 3]{ns},
\cite{scholten} concerning the hyperelliptic locus
$\mm^{\mbox{\scriptsize $\op{h}$}}_3$ we get a complete picture
for $\mm_3$:

\begin{teor}\label{mass}
The number of $\fq$-rational points of the strata by the Newton
polygon of $\mm^{\mbox{\scriptsize $\op{nh}$}}_3$,
$\mm^{\mbox{\scriptsize $\op{h}$}}_3$ and $\mm_3$ are given in the
following table:

\renewcommand\arraystretch{1.4}
\centerline{\begin{tabular}{c|ccccc}
&ordinary& 2-rank two&2-rank one&type $1/3$&supersingular\\
\hline $\mm_3^{\textrm{nh}}$& $q^6-q^5+1$&$q^5-q^4$&$q^4-q^3$&$q^3-q^2$&$q^2$\\
$\mm_3^{\textrm{h}}$&$q^5-q^4$&$q^4-2q^3+q^2$&$2(q^3-q^2)$&$q^2$&0\\
$\mm_3$& $q^6-q^4+1$&$q^5-2q^3+q^2$&$q^4+q^3-2q^2$&$q^3$&$q^2$
\end{tabular}}
\renewcommand\arraystretch{1.}
\end{teor}\e

In particular, $|\mm_3(\fq)|=q^6+q^5+1$. This result for the total
number of rational points of the moduli space was conjectured by
A. Granville and B. Brock \cite{bg} and proved by J. Bergstr\"om
\cite{be}.\e

\ni{\bf Acknowledgements. }This paper was written while the
authors enjoyed a research stay at the
 Institut f\"ur Experimentelle Mathematik in Essen. The authors
 wish to express their gratitude to the IEM
 staff and the colleagues of the
Zahlentheorie Gruppe for their warm hospitality. Also, it is a
pleasure to thank Yann Sepulcre and Xavier Xarles for some helpful
conversations. Finally, we are indebted to Bradley Brock for
calling our attention on the papers \cite{bg}, \cite{be} and for
pointing out a mistake in a previous version of the paper.

\e

{\small
\begin{tabular}{ll}
\begin{tabular}{l}
Enric Nart\\
Departament de Matem\`atiques\\ Univ. Aut\`onoma de Barcelona\\
08193 Bellaterra, Barcelona, Spain\\
{\tt nart@mat.uab.es}
\end{tabular}&
\begin{tabular}{l}
Christophe Ritzenthaler\\
Institut f\"ur Experimentelle Mathematik\\ Ellernstr. 29\\
45326 Essen, Germany\\ {\tt ritzenth@math.jussieu.fr}
\end{tabular}\end{tabular}}

\end{document}